\documentclass[journal=jacsat,manuscript=article]{achemso}
\usepackage[T1]{fontenc} 
\usepackage{amsmath}
\usepackage{amssymb}
\usepackage{amsthm}
\usepackage{graphicx}
\usepackage{mathrsfs}
\usepackage{caption}
\usepackage{subcaption}
\usepackage{overpic}
\usepackage[all]{xy}
\usepackage{tikz}
\usetikzlibrary{arrows.meta}
\usetikzlibrary{backgrounds}
\usepackage[hidelinks,breaklinks]{hyperref}

\DeclareMathOperator{\grad}{grad}
\DeclareMathOperator{\Hess}{Hess}

\DeclareMathOperator{\sgn}{sgn}

\newcommand{\R}{\mathbb{R}}
\renewcommand{\S}{\mathbb{S}}
\newcommand{\X}{\mathfrak{X}}
\renewcommand{\d}{\mathrm{d}}
\newcommand{\e}{\mathrm{e}}
\newcommand{\D}{\mathscr{D}}
\renewcommand{\epsilon}{\varepsilon}

\theoremstyle{plain}

\theoremstyle{definition}

\newtheorem{example}{Example}
\theoremstyle{remark}
\newtheorem{remark}{Remark}

\author{Juan M. Bello-Rivas}
\email{jmbr@jhu.edu}
\affiliation[JHU]{Department of Chemical and Biomolecular Engineering, Whiting School of Engineering, Johns Hopkins University, 3400 North Charles Street, Baltimore, 21218, MD, USA}
\author{Anastasia Georgiou}
\affiliation[JHU]{Department of Chemical and Biomolecular Engineering, Whiting School of Engineering, Johns Hopkins University, 3400 North Charles Street, Baltimore, 21218, MD, USA}
\author{Hannes Vandecasteele}
\affiliation[KULeuven]{Department of Computer Science, KU Leuven, Celestijnenlaan 200A, 3001 Leuven}
\author{Ioannis G. Kevrekidis}
\email{yannisk@jhu.edu}
\affiliation[JHU]{Department of Chemical and Biomolecular Engineering, Whiting School of Engineering, Johns Hopkins University, 3400 North Charles Street, Baltimore, 21218, MD, USA}
\alsoaffiliation{Departments of Applied Mathematics and Statistics, Johns Hopkins University, 3400 North Charles Street, Baltimore, 21218, MD, USA}

\title[GAD on manifolds defined by point-clouds]{Gentlest ascent dynamics on manifolds defined by adaptively sampled point-clouds}

\begin{document}

\begin{tocentry}
  \begin{center}
    \includegraphics[width=0.54\columnwidth]{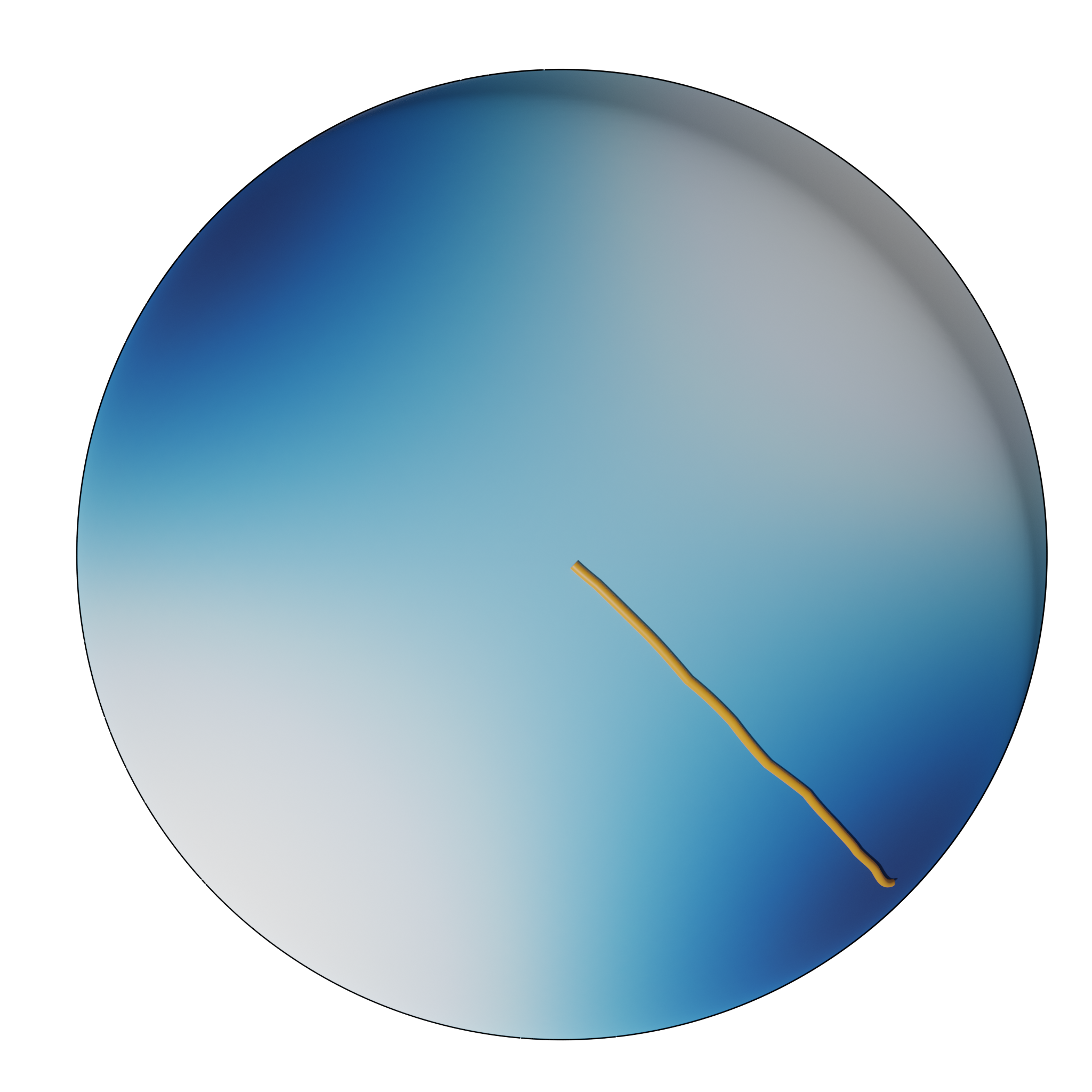}
  \end{center}
\end{tocentry}

\begin{abstract}
  Finding saddle points of dynamical systems is an important problem in practical applications such as the study of rare events of molecular systems.
  Gentlest ascent dynamics (GAD)~\cite{e2011} is one of a number of algorithms in existence that attempt to find saddle points in dynamical systems.
  It works by deriving a new dynamical system in which saddle points of the original system become stable equilibria.
  GAD has been recently generalized to the study of dynamical systems on manifolds (differential algebraic equations) described by equality constraints~\cite{yin2022a} and given in an extrinsic formulation.
  In this paper, we present an extension of GAD to manifolds defined by point-clouds, formulated using the intrinsic viewpoint.
  These point-clouds are adaptively sampled during an iterative process that drives the system from the initial conformation (typically in the neighborhood of a stable equilibrium) to a saddle point.
  Our method requires the reactant (initial conformation), does not require the explicit constraint equations to be specified, and is purely data-driven.
\end{abstract}

\section{Introduction}

The problem of finding saddle points of dynamical systems has one of its most notable applications in the search for transition states of chemical systems described at the atomistic level, since saddle points coincide with transition states~\cite{hanggi1990} at the zero temperature limit.
While finding local stable equilibria (sinks) is a relatively straightforward matter, finding saddle points is a more complicated endeavor for which a number of algorithms have been presented in the literature~\cite{olsen2004}.

Saddle point search methods can be classified according to whether they require one single input state (usually the reactant, located at the minimum of a free energy well) or two states (reactant \emph{and} product).
The Gentlest Ascent Dynamics (GAD)~\cite{e2011} belongs to the class of methods requiring a single reactant state as input and its applicability has been demonstrated with atomistic chemical systems~\cite{samanta2014,quapp2014}.
GAD can be regarded as a variant of the dimer method that is formulated as a continuous dynamical system whose integral curves with initial condition at the reactant state can lead to saddle points.
Variants of GAD such as high-index saddle dynamics (HiSD) have been the subject of recent research efforts and applications~\cite{gu2017,gu2018,yin2019,yin2022a,yin2022b,gu2022,luo2022,zhang2022}.

While many search schemes attempt to find an optimal path between reactant and product (or between reactant state and transition state), it is interesting that there exist continuous curves joining the desired states in a variety of ways: following isoclines~\cite{quapp1998,quapp1998a}, gradient extremals~\cite{basilevsky1982,hoffman1986}, and the GAD studied in this paper.
In most cases the study of these curves has been carried out in the Euclidean setting with some exceptions on the manifold of internal coordinates~\cite{quapp2004,quapp2014} and on manifolds defined by the zeros of smooth maps~\cite{yin2022a} (in these cases, the algorithms are formulated extrinsically on the ambient space).
Our contribution is formulated intrinsically and is valid on arbitrary manifolds, not necessarily explicitly defined by an atlas or by the zeros of maps.
Importantly, algorithms like the one presented here or our previous work~\cite{bello-rivas2022} do not rely on \emph{a priori} knowledge of good collective coordinates, but rather use manifold learning to find them on the fly.
In our method, there is a feedback loop of data collection that drives progress towards a saddle point.

In this paper, we study an application of the GAD to manifolds defined by point-clouds.
The manifold does not need to be characterized in advance either by the zeros of a smooth function or by an atlas, and it is only assumed that the user is capable of sampling the vicinity of arbitrary points on the manifold (\emph{e.g.,} umbrella sampling based on reduced local coordinates). 
The algorithm uses dimensionality reduction (namely, diffusion map coordinates~\cite{coifman2006}) to define a dynamical system intrinsically on the reduced coordinates that can lead to a saddle point.
Since the saddle point in general is not expected to lie in the vicinity of the reactant, our algorithm works by iteratively sampling the manifold on the fly, resolving the path on the local chart, and repeatedly switching charts until convergence.
Our approach shares algorithmic elements with our previous work~\cite{bello-rivas2022} which, however, instead of GAD dynamics on manifolds, was following isoclines on manifolds.


\section{Gentlest Ascent Dynamics and Idealized Saddle Dynamics}
\label{sec:gad}

Gentlest Ascent Dynamics (GAD)~\cite{e2011} is an algorithm for finding saddle points of dynamical systems.
We propose an extension of GAD to manifolds defined by point-clouds that finds saddle points combining nonlinear dimensionality reduction and adaptive sampling.

Let $U \colon \R^n \to \R$ be a smooth potential energy function and consider the associated gradient vector field $X$ in $\R^n$ given by $\dot{x} = X(x)$, where $X = -D U$ with $D$ denoting the gradient (or, equivalently, the Jacobian matrix).
We restrict ourselves here, for the sake of simplicity, to the case of gradient systems and index-1 saddle points.
The GAD algorithm consists of integrating the equations of motion of the related dynamical system $\hat{X}$ on an extended phase space $\R^{2n}$ given by
\begin{equation}
  \left\{
    \begin{aligned}
      &\dot{x} = -H(v) D U(x), \\
      &\dot{v} = -D^2 U(x) v + r(x, v) v,
    \end{aligned}
  \right.
  \label{eq:gad-orig}
\end{equation}
where $x, v \in \R^n$, $H(v) w = w - (2 \, v \! \cdot \! w) v$ is the Householder reflection~\cite{golub2013} of $w \in \R^n$ across the hyperplane $\langle v \rangle^\perp = \{ z \in \R^n \mid v \cdot z = 0 \}$ and $r(x, v) = \| v \|^{-2} \, v \cdot D^2 U(x) v$ is the Rayleigh quotient of the Hessian matrix $D^2 U(x)$ corresponding to the vector $v \in \R^n$.

\begin{remark}
  The right hand side of the ordinary differential equation for $\dot{v}$ is the term
  \begin{equation*}
    -D^2 U(x) v + r(x, v) v.
  \end{equation*}
  If $v$ is an eigenvector of $D^2 U(x)$ with eigenvalue $\lambda$, then $r(x, v) = \lambda$.
  Now note that the constrained optimization problem consisting of finding the extrema of $Z(x) = \frac{1}{2} \| -D U(x) \|^2$ along the level sets $C = \{ x \in \R^n \mid U(x) = c \}$ is precisely given by the Lagrange equation
  \begin{equation*}
    D^2 U(x) D U(x) - \lambda D U(x) = 0,
  \end{equation*}
  which says that the extrema of the magnitude of the gradient, $Z$, along $C$ are attained wherever the gradient field $X = -D U$ happens to be an eigenvector of the Hessian $D^2 U$.
\end{remark}

\begin{figure}[ht]
  \centering
  \includegraphics{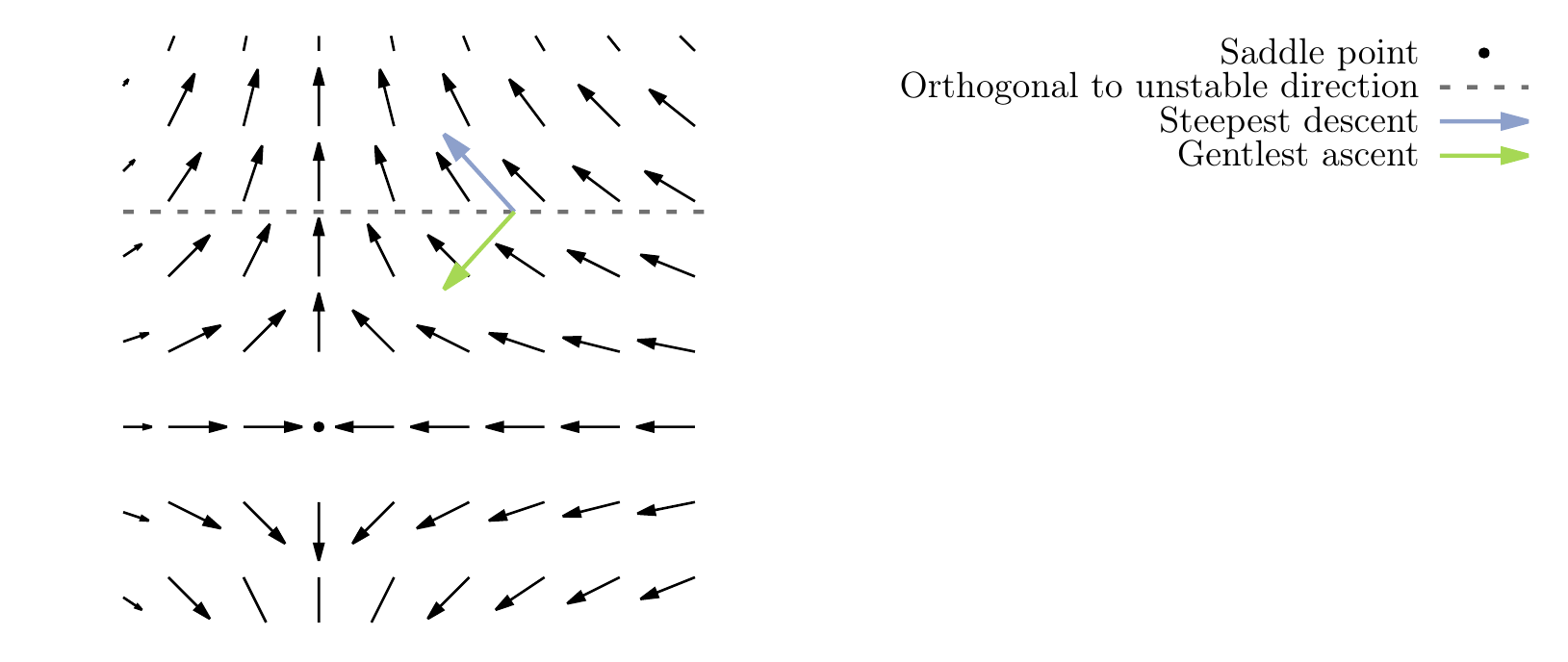}
  \caption{Vector field corresponding to the first equation of GAD on the potential $U(x^1, x^2) = (x^1)^2 - (x^2)^2$ with $v = (0, 1) \in \R^2$ being the unstable direction. Note how the steepest descent direction is reflected across the orthogonal complement $\langle v \rangle^\perp$ to obtain $\dot{x}$.}
  \label{fig:gentlest-ascent}
\end{figure}

The intuition behind GAD is that the force $-D U$ can be written in a basis determined by the eigenvectors of the Hessian $D^2 U$, which gives us the stable and unstable directions in the vicinity of an index-1 saddle point.
The second differential equation in~\eqref{eq:gad-orig} acts as a (continuous) eigensolver yielding the unstable direction.
Then, we can flip the sign of the component of the force corresponding to the unstable direction to make the resulting vector point towards the saddle point (see Figure~\ref{fig:gentlest-ascent}).

\subsection{Idealized saddle-point dynamics}
\label{sec:isd}

We consider a variant of the GAD algorithm named Idealized Saddle Dynamics (ISD)~\cite{levitt2017}, given by
\begin{equation*}
  \dot{x} = -H(v) D U(x),
\end{equation*}
where $v \in \R^n$ is an eigenvector of $D^2 U(x)$ corresponding to the smallest eigenvalue $\lambda$.

We now leave the Euclidean setting behind and study the Riemannian setting (see the appendix for a summary of the topic and the notation).
Let $M$ be a $d$-dimensional smooth manifold with Riemannian metric $g$ and let $U \in C^\infty(M)$ be a potential energy function.
We consider the gradient field $X = -\grad U \in \X(M)$ and we write the following equation for the ISD vector field on $M$,
\begin{equation}
  \label{eq:isd}
  \hat{X}
  =
  X - 2 g(V, X) V
  \in
  \X(M),
\end{equation}
where $V$ is the vector field on $M$ defined by choosing the eigenvector (normalized so that $g(V, V) = 1$) corresponding to the smallest eigenvalue of the Hessian matrix.
However, in this case the Hessian must be defined in terms of the covariant derivative on $(M, g)$ with respect to the Levi-Civita connection~\cite{docarmo1992,absil2008}.
To be precise, given two vector fields $S, T \in \X(M)$, we have
\begin{equation*}
  \nabla^2 U(S, T)
  =
  \nabla_T (\nabla_S U) - \nabla_{\nabla_T S} U,
\end{equation*}
which is a tensor of type $(0, 2)$, where $\nabla$ denotes the covariant derivative induced by the Levi-Civita connection (the notation $\nabla_T S$ represents the covariant derivative of $S$ in the direction of $T$).
We apply the sharp ($\sharp$) isomorphism (raising or lowering indices) to turn $\nabla^2 U$ into a $(1, 1)$-tensor.
Therefore,
\begin{equation*}
  (\Hess U) T
  =
  \nabla_T \grad U.
\end{equation*}
The eigenvector corresponding to the lowest eigenvalue of $\Hess U$ at each point $p \in M$ induces a vector field $V$ whose integral curves are curves that may join an initial point with a saddle point.

The ISD formulation of GAD is particularly amenable to be coupled with dimensionality reduction approaches because the resulting eigenproblem is often of much lower dimensionality than the ambient space in which the original dynamical system is defined.

\begin{example}[Exact solution of a model system]
  \label{ex:isd-simple-stereographic}
  Let us compute a concrete case of ISD, first with known exact formulas and later on with approximations using diffusion maps and Gaussian processes.
  Consider the sphere
  \begin{equation*}
    \S^2
    =
    \{ (x^1, x^2, x^3) \in \R^3 \mid (x^1)^2 + (x^2)^2 + (x^3)^2 = 1 \}.
  \end{equation*}
  and the stereographic projection from the North pole onto the tangent plane at the South pole.
  The system of coordinates is given by
  \begin{equation*}
    \phi(x^1, x^2, x^3) = (\frac{x^1}{1 - x^3}, \frac{x^2}{1 - x^3}) \in \R^2,
  \end{equation*}
  Let $(u^1, u^2) \in \R^2$.
  The corresponding parameterization, $\psi = \phi^{-1}$, is given by
  \begin{equation*}
    \psi(u^1, u^2)
    =
    \frac{1}{1 + (u^1)^2 + (u^2)^2} (2 u^1, 2 u^2, (u^1)^2 + (u^2)^2 - 1).
  \end{equation*}
  The pullback of the Euclidean metric $h = \d x^1 \otimes \d x^1 + \d x^2 \otimes \d x^2 + \d x^3 \otimes \d x^3$ by $\psi$ gives us the metric $g$
  \begin{equation*}
    g
    =
    \psi^\star h
    =
    \frac{4 \d u^1 \otimes \d u^1 + 4 \d u^2 \otimes \d u^2}{\left( 1 + (u^1)^2 + (u^2)^2 \right)^2}
  \end{equation*}
  The non-redundant Christoffel symbols $\Gamma_{ij}^k$ that characterize the Levi-Civita connection $\nabla$ are
  \begin{equation*}
    \Gamma_{1, 1}^{1}
    =
    \frac{-2 u^1}{1 + (u^1)^2 + (u^2)^2},
    \quad
    \Gamma_{1, 2}^{1}
    =
    \frac{-2 u^2}{1 + (u^1)^2 + (u^2)^2},
  \end{equation*}
  \begin{equation*}
    \Gamma_{2, 2}^{1} = -\Gamma_{1, 1}^{1},
    \quad
    \Gamma_{1, 1}^{2} = -\Gamma_{1, 2}^{1},
    \quad
    \Gamma_{1, 2}^{2} = \Gamma_{1, 1}^{1},
    \quad \text{and} \quad
    \Gamma_{2, 2}^{2} = \Gamma_{1, 2}^{1}.
  \end{equation*}
  The energy $E(x^1, x^2, x^3) = x^1 x^2 x^3$, constrained on $\S^2 \subset \R^3$ is transformed to $U = \psi^\star E$ on $\S^2$,
  \begin{equation*}
    U(u^1, u^2)
    =
    \frac{4 u^1 u^2 \left( (u^1)^2 + (u^2)^2 - 1 \right)}{\left( (u^1)^2 + (u^2)^2 + 1 \right)^3},
  \end{equation*}
  The force on $\S^2$ is the negative of the gradient
  \begin{multline*}
    \grad U
    =
    \left( {{(u^{2})^5-2\,(u^{1})^2\,(u^{2})^3+\left(-3\,(u^{1})^4+8\,(u^{1})^2-1\right)\,{u^{2}}}\over{(u^{2})^4+\left(2\,(u^{1})^2+2\right)\,(u^{2})^2+(u^{1})^4+2\,(u^{1})^2+1}}  \right) \, \frac{\partial}{\partial {u^{1}}} \\
    +
    \left( -{{3\,{u^{1}}\,(u^{2})^4+\left(2\,(u^{1})^3-8\,{u^{1}}\right)\,(u^{2})^2-(u^{1})^5+{u^{1}}}\over{(u^{2})^4+\left(2\,(u^{1})^2+2\right)\,(u^{2})^2+(u^{1})^4+2\,(u^{1})^2+1}}  \right) \, \frac{\partial}{\partial {u^{2}}}
  \end{multline*}
  The potential energy and the gradient field are shown in Figure~\ref{fig:isd-simple-sphere-a}.

  \begin{figure}[ht]
    \centering
    \begin{subfigure}[t]{0.75\columnwidth}
      \includegraphics[trim=1cm 0.5cm 1.25cm 1cm,clip]{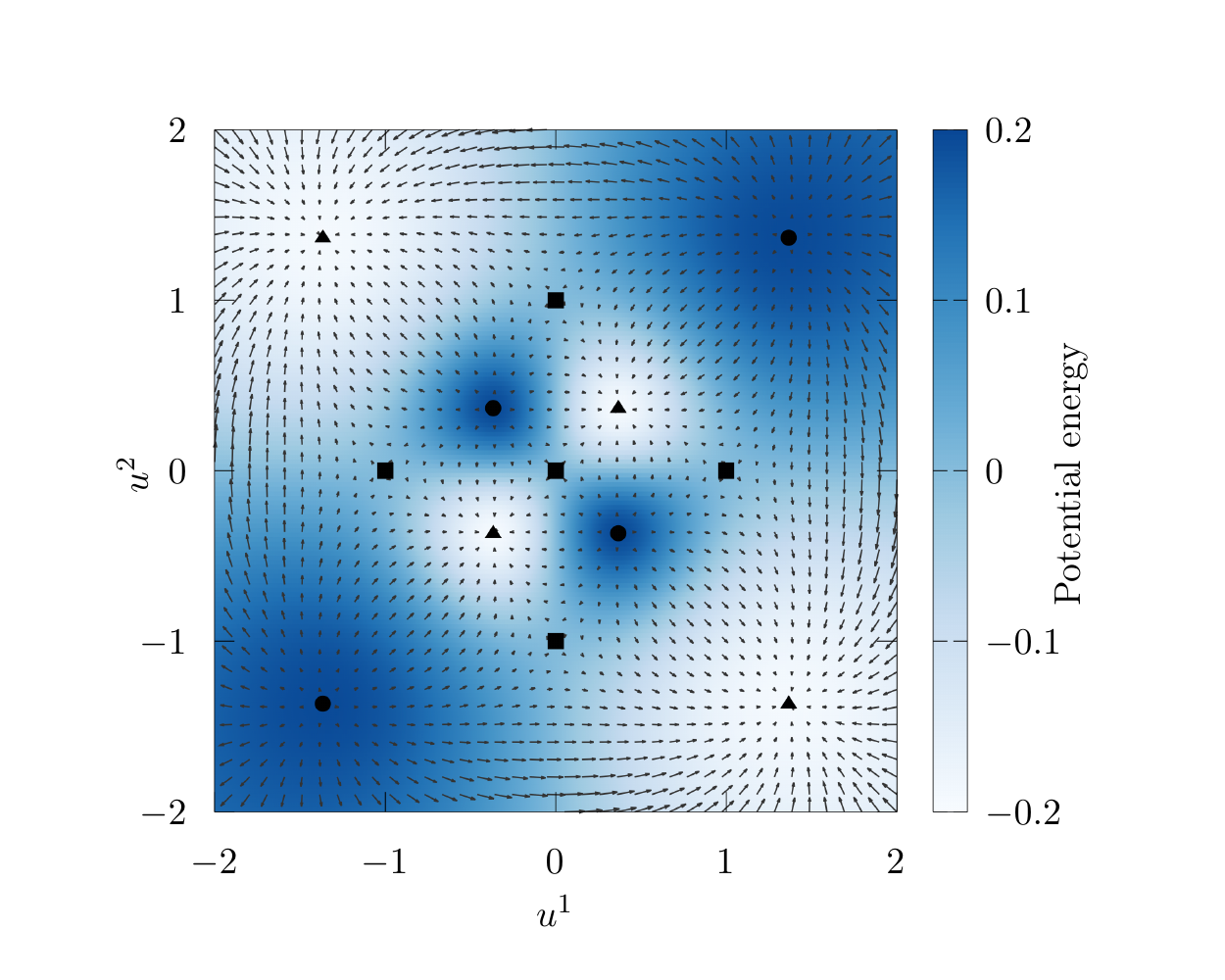}
      \caption{Vector field corresponding to the potential energy function considered in the text.}
      \label{fig:isd-simple-sphere-a}
    \end{subfigure}
    \begin{subfigure}[t]{0.75\columnwidth}
      \includegraphics[trim=1cm 0.5cm 1.25cm 1cm,clip]{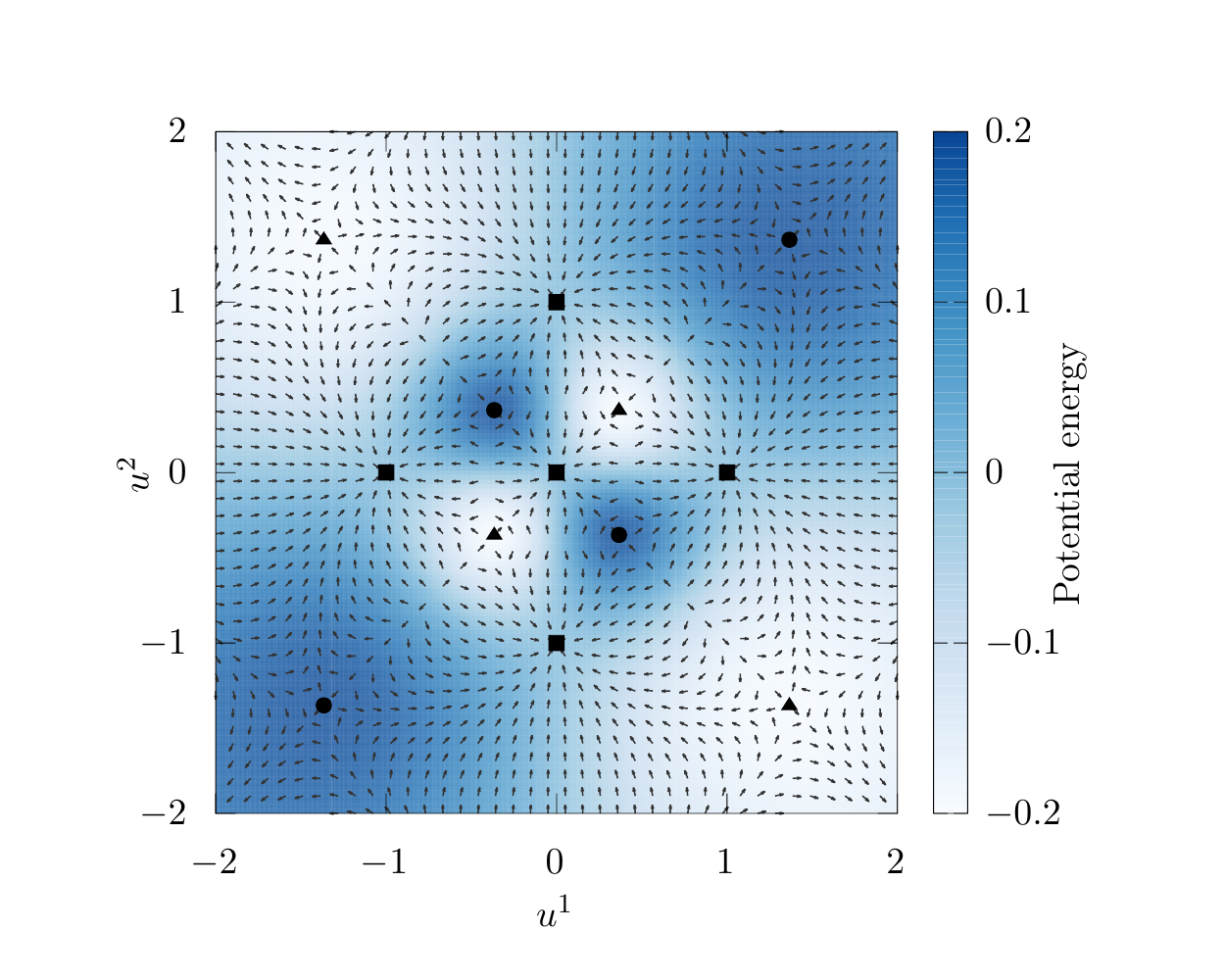}
      \caption{Idealized saddle dynamics vector field corresponding to the original vector field. Notice how saddle points in the original vector field become sinks (stable equilibria) of the ISD vector field.}
      \label{fig:isd-simple-sphere-b}
    \end{subfigure}
    \caption{Vector fields on the sphere, seen in the stereographic projection from the North pole onto the tangent plane to the South pole. The equilibria of the original vector field appear represented as $\bullet$ for sinks, $\blacksquare$ for saddle points, and $\blacktriangle$ for sources.}
    \label{fig:isd-simple-sphere}
  \end{figure}

  The Hessian is then given by the $(1, 1)$ tensor
  \begin{equation*}
    \nabla \grad U
    =
    A
    \frac{\partial}{\partial {u^{1}}} \otimes \d u^1
    +
    B
    \frac{\partial}{\partial {u^{2}}} \otimes \d u^1
    +
    B
    \frac{\partial}{\partial {u^{1}}} \otimes \d u^2
    +
    D
    \frac{\partial}{\partial {u^{2}}} \otimes \d u^2,
  \end{equation*}
  where
  {\footnotesize
  \begin{align*}
    A
    &=
-{{4\,{u^{1}}\,{u^{2}}\,\left((u^{2})^4+(u^{2})^2-(u^{1})^4+11\,(u^{1})^2-6\right)}\over{\left((u^{2})^2+(u^{1})^2+1\right)^3}} \\
    B
    &=
      - {{(u^{2})^6-5\,(u^{1})^2\,(u^{2})^4-5\,(u^{2})^4-5\,(u^{1})^4\,(u^{2})^2+30\,(u^{1})^2\,(u^{2})^2-5\,(u^{2})^2+(u^{1})^6-5\,(u^{1})^4-5\,(u^{1})^2+1}\over{\left((u^{2})^2+(u^{1})^2+1\right)^3}} \\
    D
    &=
{{4\,{u^{1}}\,{u^{2}}\,\left((u^{2})^4-11\,(u^{2})^2-(u^{1})^4-(u^{1})^2+6\right)}\over{\left((u^{2})^2+(u^{1})^2+1\right)^3}}.
  \end{align*}
  }
  The smallest eigenvector of $\nabla \grad U$ determines the vector field $V \in \X(\S^2)$ that is used in the formulation of the ISD vector field~\eqref{eq:isd} and is shown in Figure~\ref{fig:isd-simple-sphere-b}.
\end{example}

Example~\ref{ex:isd-simple-stereographic} required the knowledge of a particular system of coordinates (namely, the stereographic projection from the North pole to the tangent space at the South pole) mapping three-dimensional points on the sphere to two-dimensional coordinates; here this allows us to work with a single chart.
In some settings the system of coordinates of the underlying manifold is either unknown or difficult to obtain.
It is possible under those circumstances to replicate the steps in Example~\ref{ex:isd-simple-stereographic} in the absence of a system of coordinates given as a closed-form expression by resorting to manifold learning / dimensionality reduction techniques.
In our case, as we shall discuss next, we use diffusion maps on points sampled from a local neighborhood of the manifold to extract a suitable system of coordinates.
Fitting a Gaussian process to the diffusion map coordinates of the point-cloud yields a local system of coordinates $\phi$ that can be evaluated at arbitrary points (not necessarily those in the sampled point-cloud).
Once we have a system of coordinates, we can again proceed to estimate the Riemannian metric as well as the Levi-Civita connection, and compute the flow of the ISD vector field $\hat{X}$ given in~\eqref{eq:isd} to find saddle points.

Revisiting Example~\ref{ex:isd-simple-stereographic} and approaching it with the aforementioned procedure allows us to obtain trajectories that lead to saddle points, as shown in Figures~\ref{fig:isd-diffusion-maps-1} to~\ref{fig:isd-diffusion-maps-10}.

The choice of diffusion map coordinates for dimensionality reduction and Gaussian processes for non-linear regression is inessential and both methods could be replaced by, for instance, neural networks (\emph{e.g.,} a variational autoencoder and a graph neural network, respectively).
An important aspect of our approach is the fact that it operates on local neighborhoods of points, as opposed to using data from longtime simulations, and only then applying dimensionality reduction techniques to extract global collective variables.
From a geometric viewpoint, this is motivated by the fact that a sufficiently small neighborhood of a manifold can always be transformed onto a subset of its tangent space by a smooth invertible map.
From a physical viewpoint, different sets of collective variables govern different stages of a reaction (\emph{e.g.,} the distance between a ligand and a receptor is important when both molecules are far away whereas their relative orientation might be more important to ligand-binding, when both molecules become close).

One interesting aspect of diffusion maps is its relation to the infinitesimal generator of a diffusion process~\cite{nadler2006}, which in turn has connections to the committor function, which is an optimal reaction coordinate~\cite{peters2016,berezhkovskii2019,roux2022}.
We conclude this discussion on collective variables by pointing out that \emph{a priori} knowledge of good reaction coordinates (some recent examples can be found in~\cite{wu2022,manuchehrfar2022}) can often be put in one-to-one correspondence with diffusion map coordinates~\cite{chiavazzo2017,fujisaki2018}.

\subsection{Mean force}

In computational statistical mechanics we are often interested in dynamics on Riemannian ma\-nifolds endowed with a probability distribution.
For instance, for simulations at constant temperature, we use the Boltzmann distribution with probability density proportional to $\exp\{-\beta U(x)\}$, where $\beta > 0$ is the inverse temperature.
The relevant vector field at $u$ in the local chart is the \emph{mean force}~\cite{fixman1974,carter1989,denotter2000}, given by
\begin{equation*}
  X_u
  =
  -
  \langle
  \grad
  \left(
    U \circ \psi
    -
    \tfrac{1}{2}
    \beta^{-1}
    \log \det (D \psi^\top D \psi)
  \right)
  \rangle_u,
\end{equation*}
where $\langle \cdot \rangle_u$ denotes the ensemble average with respect to the Boltzmann distribution, conditioned on $\{ x \in \R^n \mid \phi(x) = u \}$.
There are a number of numerical methods that estimate the mean force as a step in their calculations.
Adaptive biasing force~\cite{darve2008} and the string method~\cite{maragliano2006} are examples of those.


\section{Algorithm}
\label{sec:algorithm}

Consider the manifold $M \subset \R^n$ and the gradient dynamical system $X \in \X(M)$.
We begin by drawing a total of $N$ samples from the manifold $M$ in the neighborhood of an initial point $p \in M$.
This can be done in a variety of ways depending on the application.
If $M$ is the inertial manifold of a dynamical system $X$, then a reasonable way to approach this problem is to generate $N$ distinct perturbations $\{ p_{(i)} \in \R^n \mid i = 1, \dotsc, N \}$ of $p$ and propagate them according to the flow of $X$ during a (short) time horizon $\tau > 0$.
Doing so, we obtain a data set $\D = \{ q_{(i)} = \exp_\tau p_{(i)} \in M \mid i = 1, \dotsc, N \}$ approximately on the manifold. 
Alternatively, one may numerically solve a stochastic differential equation such as the Brownian dynamics equation,
\begin{equation}
  \label{eq:sde}
  \d q_t
  =
  X(q_t) \, \d t
  +
  \sigma \, \d B_t,
\end{equation}
where the drift is the vector field $X$, $\sigma > 0$ is a constant, and $B_t$ is a standard $n$-dimensional Brownian motion.
Solving~\eqref{eq:sde} (possibly with an added RMSD-based restraint around the initial conformation) up to a certain time $\tau > 0$ and extracting an uncorrelated subset of the states at different time steps yields a data set $\D = \{ q_{(i)} = q_{t_i} \mid i = 1, \dotsc, N, 0 \le t_1 \le \cdots \le t_N \le \tau \}$.

Next, we apply a dimensionality reduction algorithm on the data set $\D$ to obtain a set of reduced coordinates $\phi$.
In our case, we use diffusion maps~\cite{coifman2006,nadler2006} to obtain a set of vectors $\phi_{(i)} \in \R^d$ with $d \le n$ but other methods, such as local tangent space alignment~\cite{zhang2004}, may be used as well.
It is important to note that our dimensionality reduction method is applied to a \emph{local neighborhood} of an initial point and, therefore, it is expected to yield a reasonable approximation to a chart \emph{on that neighborhood}.

We fit a Gaussian process regressor $\phi$ to the pairs of points $(q_{(i)}, \phi_{(i)}) \in \R^n \times \R^d$ to obtain a smooth map $\phi \colon M \subset \R^n \to \R^d$ that will act as a system of coordinates (in particular, $\phi(q_{(i)}) \approx \phi_{(i)}$).
Proceeding in an analogous fashion, we compute the inverse mapping $\psi = \phi^{-1}$.

Note that one possible way of estimating the dimensionality is by computing the average of the approximate rank of the Jacobian matrix of $\phi$ at (a subset of) the data points and retaining the components that yield a local chart.

\begin{remark}
  We can reduce the computational expense of the Gaussian process regression by reusing the kernel matrix with entries $\e^{-\| q_{(i)} - q_{(j)} \|^2 / 2 \epsilon}$ (for some $\epsilon > 0$) calculated during the computation of the diffusion map coordinates as the covariance matrix for the Gaussian process (assuming that it is formulated using the squared exponential kernel).
\end{remark}

\begin{remark}
  It is not always possible to obtain a suitable Gaussian process regressor for $\psi = \phi^{-1}$.
  An alternative is to add an Ornstein-Uhlenbeck process to the stochastic differential equation~\eqref{eq:sde} in order to obtain
  \begin{equation*}
    \d q_t
    =
    \left(
      X(q_t) - \kappa (\phi(q_t) - \phi_0)
    \right) \, \d t
    +
    \sigma \, \d B_t,
  \end{equation*}
  where $\kappa > 0$ is a hyperparameter and $\phi_0$ is a prescribed point not necessarily in $\phi(\D)$.
  Computing the ensemble average $\langle q_t \rangle$ of the solution to the above equation yields a point $q_{(0)}$ such that $\phi(q_{(0)}) \approx \phi_0$ or, equivalently, $q_{(0)} = \phi^{-1}(\phi_{(0)})$.
  This works because the new term added to the drift nudges the system towards a point in ambient space such that its image by $\phi$ is the prescribed point $\phi_0$.
\end{remark}

We consider the values $X_{q_{(i)}} \in T M$ of the vector field $X$ at the points in the data set $\D$ and map them via the system of coordinates in order to obtain the vector field $\phi_\star X_{q_{(i)}} = D \phi(q_{(i)}) X_{q_{(i)}}$ in the new coordinates, where $D \phi$ is the Jacobian matrix of $\phi$ (note that the Jacobian-vector product can be computed either as a closed-form formula or efficiently using automatic differentiation ---\emph{e.g.,} using \texttt{jvp} in JAX~\cite{jax2018github}).
We fit another Gaussian process to the pushforward vector field $\phi_\star X$ in order to be able to evaluate it at arbitrary points in the new coordinates and we abuse notation in what follows by also referring to $\phi_\star X$ as $X$.

At this stage, we can readily compute the Riemannian metric $g$ as
\begin{equation}
  \label{eq:metric}
  g
  =
  \sum_{i, j = 1}^d g_{ij} \, \d \phi^i \otimes \d \phi^j,
\end{equation}
where $g_{ij} = D \psi^\top D \psi$ for $i, j = 1, \dotsc, d$.
The metric induces an inner product, denoted by the bracket $\langle \cdot, \cdot \rangle$, between tangent vectors such that if $S = \sum_{i = 1}^d S^i \frac{\partial}{\partial \phi^i}$ and $T = \sum_{i = 1}^d T^i \frac{\partial}{\partial \phi^i}$ are the expressions in local coordinates of two tangent vectors $S, T \in T_u M$ at a point $u \in M$, then $\langle S, T \rangle = \sum_{i, j = 1}^d g_{ij} S^i T^i$.

Using~\eqref{eq:metric}, we obtain the coefficients of the Levi-Civita connection~\cite{docarmo1992},
\begin{equation*}
  \Gamma^\ell_{jk}
  =
  \sum_{i = 1}^d g^{\ell i} \left( \frac{\partial  g_{ij}}{\partial \phi^k} + \frac{\partial  g_{ik}}{\partial \phi^j} - \frac{\partial g_{jk}}{\partial \phi^i} \right)
\end{equation*}
for $j, k, \ell \in \{ 1, \dotsc, d \}$, where $g^{ij}$ denotes the entries of the inverse of the matrix with components $g_{ij}$.
This, in turn, allows us to take the covariant derivative and the Hessian, defined by
\begin{equation*}
  \Hess U (S, T)
  =
  \langle \nabla_T \grad U, S \rangle
\end{equation*}
for arbitrary tangent vectors $S, T$.
Observe that the Hessian is defined on the local chart, not on the ambient space.
The eigenvector $V$, normalized with respect to $g$, with smallest eigenvalue of the Hessian then yields a vector field
\begin{equation*}
  \hat{X}
  =
  X
  - 2 \, V g(V, X)
\end{equation*}
such that the index-1 saddle points of $X$ become stable equilibria~\cite{e2011} of $\hat{X}$.
Consequently, an integral curve of $\hat{X}$ in the vicinity of a saddle point leads to the saddle point.

In general, in order to carry out the computation until convergence, we must frequently switch charts.
This is due to the fact that at each step, we sample a point-cloud $\D$ in a small neighborhood of a given point and the Gaussian process regressor $\phi$ yields an approximation to $X$ on the chart that cannot extrapolate far away from the sampled points.
Therefore, when we reach the confines of $\D$ (which can be determined by the density of points), we ought to map the latest point of our integral curve of $\hat{X}$ from the chart back to the ambient space using $\psi = \phi^{-1}$, as discussed earlier.
After that, we start the whole procedure again: sample a new data set $\D$, compute $\phi$, trace an integral curve of $\hat{X}$, etc.

\begin{remark}
  An alternative approach to the one presented here could consist of exploiting the fact that GAD trajectories are geodesics~\cite{bofill2016} of a Finsler metric~\cite{randers1941,chern2005} and to numerically compute said geodesics in each learned local chart in a similar manner to what we propose.
\end{remark}

The factor that impacts the algorithm the most is the deviation of the sampled points in the data set from the manifold $M$.
In other words, the farther away from $M$ our data points are, the more noisy is the estimate of $\hat{X}$, and, consequently, the harder it is to find the saddle points.

\subsection{Numerical examples}

The computations that follow were carried out using the JAX~\cite{jax2018github} and Diffrax~\cite{kidger2021} libraries, and the code to reproduce our results is available at \url{https://github.com/jmbr/gentlest_ascent_dynamics_on_manifolds}.

In our two examples below, we automatically set the bandwidth parameter (usually denoted by $\epsilon$) in diffusion maps to be equal to the squared median of the pairwise distances between points.
The regularization parameters in the Gaussian processes $\phi$ and $\psi$ are chosen so that their scores are sufficiently close to 1 on a test set.

\subsubsection{Sphere}

\begin{figure}[ht]
  \includegraphics{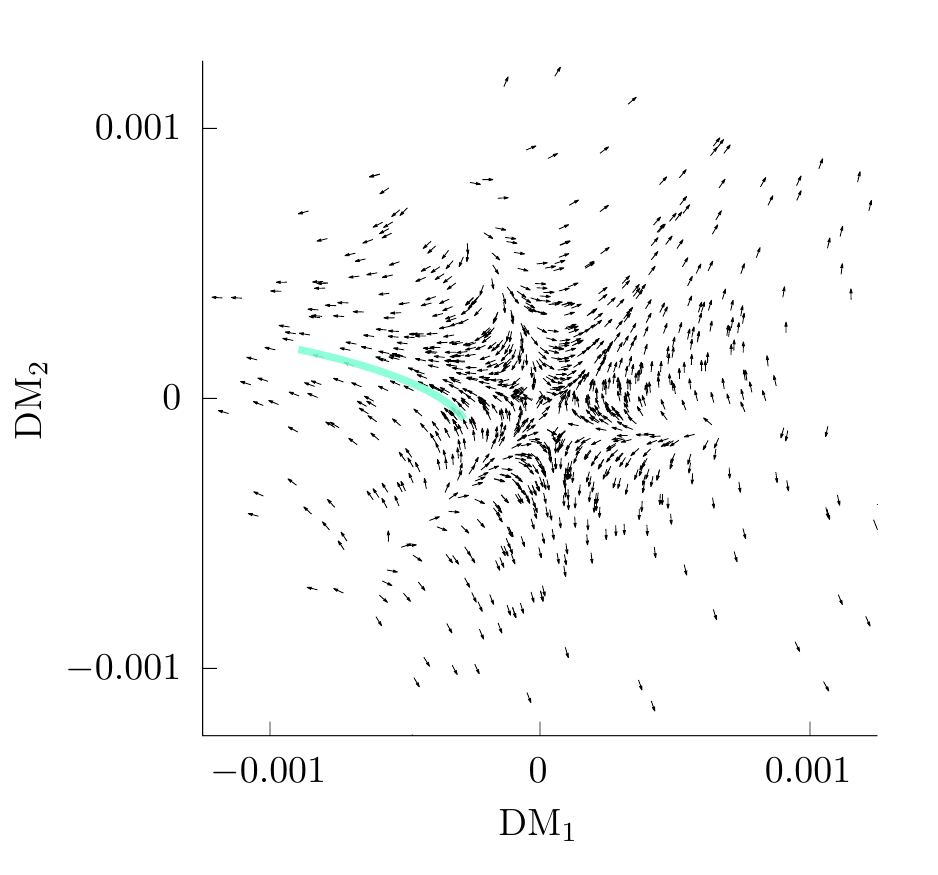}
  \caption{Iteration \#1 of ISD on the sphere by sampling point-clouds and using diffusion maps. The points and the vector field are drawn from the neighborhood of a stable equilibrium. The solid curve represents the GAD/ISD path.}
  \label{fig:isd-diffusion-maps-1}
\end{figure}

\begin{figure}[ht]
  \includegraphics{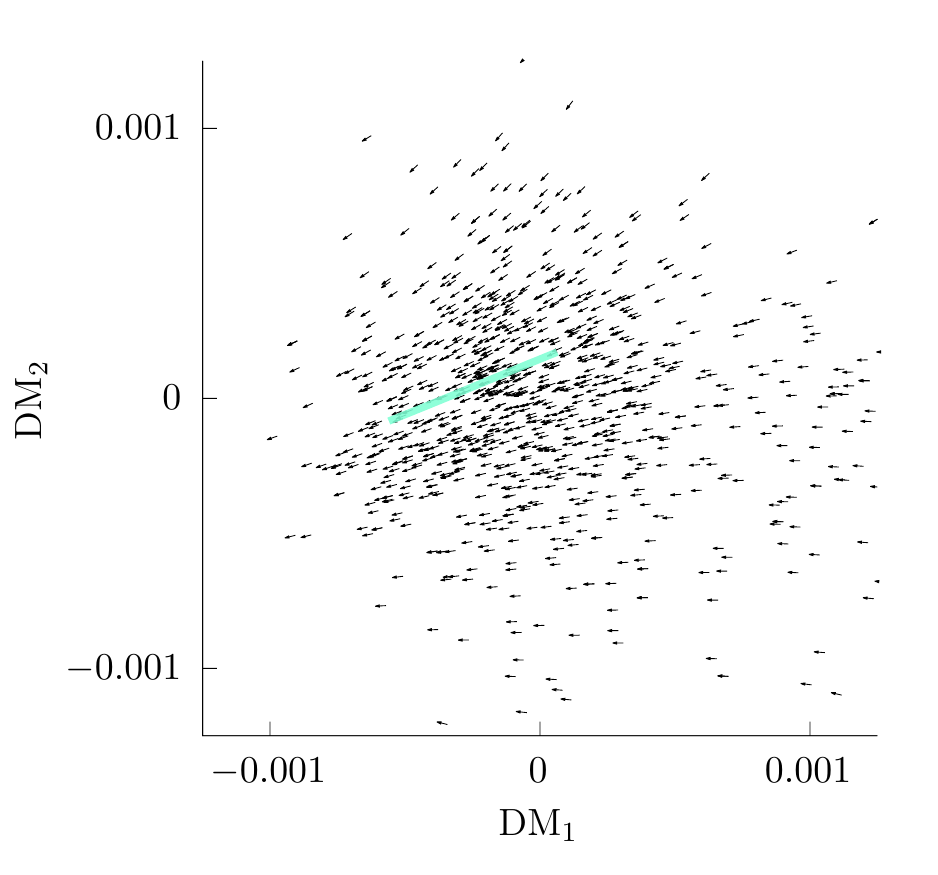}
  \caption{Iteration \#5 of ISD on the sphere by sampling point-clouds and using diffusion maps.}
  \label{fig:isd-diffusion-maps-5}
\end{figure}

\begin{figure}[ht]
  \includegraphics{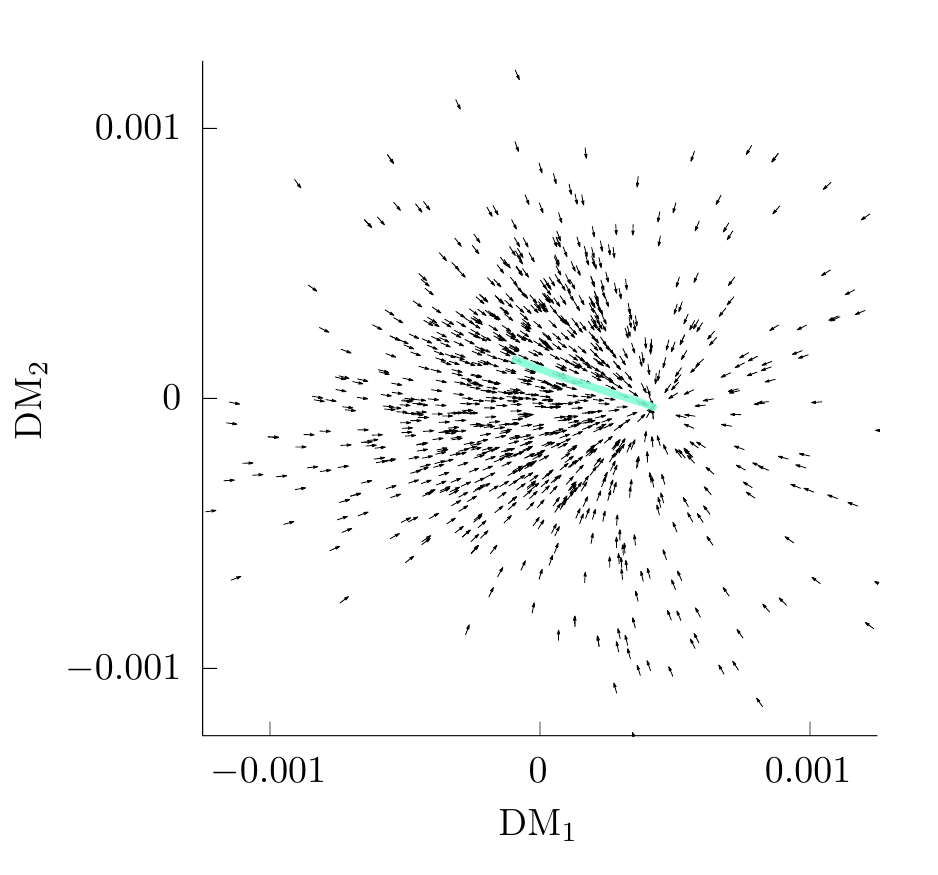}
  \caption{Iteration \#10 of ISD on the sphere by sampling point-clouds and using diffusion maps. The algorithm has reached a saddle point of the original vector field (notice how it becomes a sink for the ISD vector field).}
  \label{fig:isd-diffusion-maps-10}
\end{figure}

\begin{figure}[ht]
  \centering
  \includegraphics[width=\columnwidth]{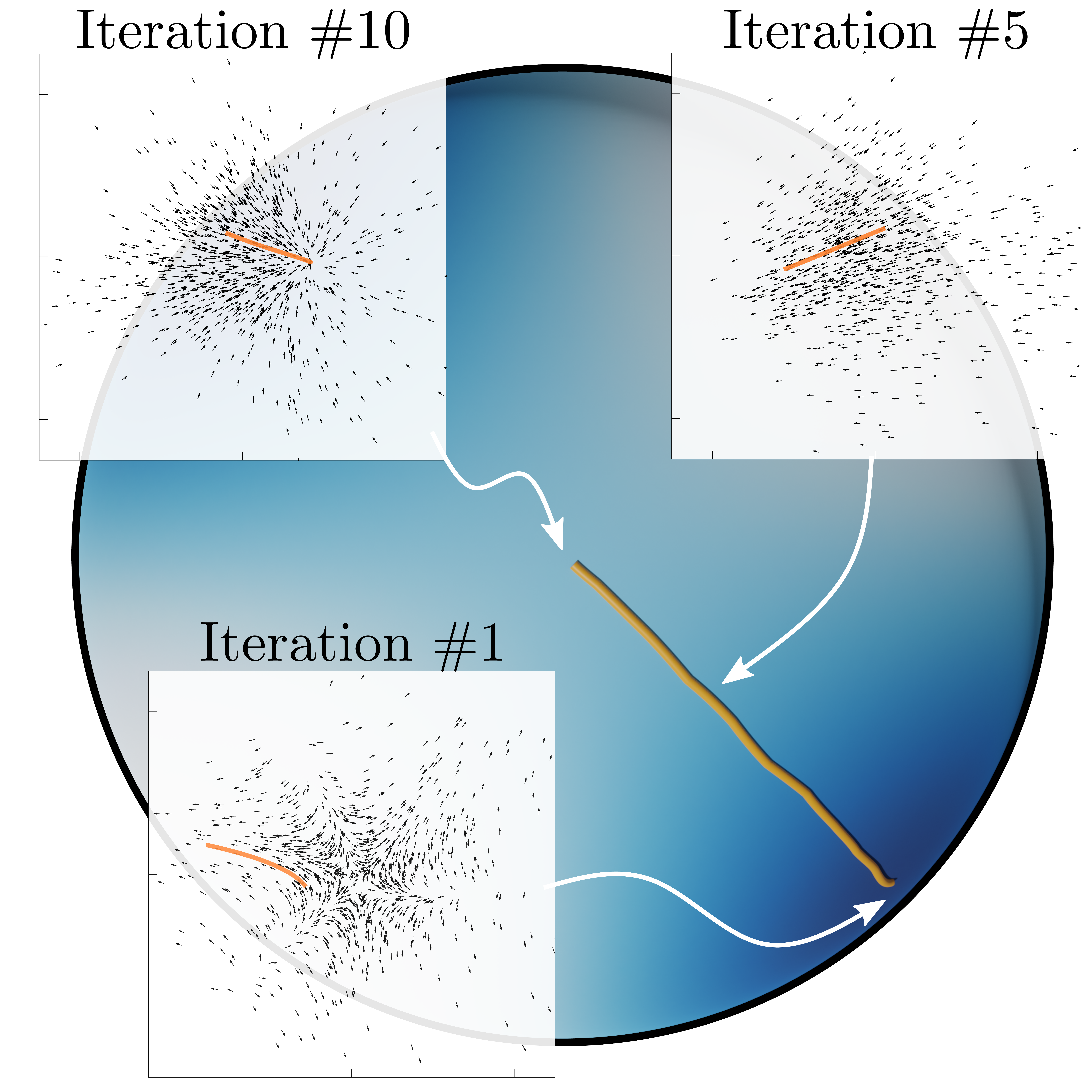}
  \caption{Resulting trajectory (in orange) of GAD/ISD on the sphere (color indicates value of the potential energy function) using diffusion maps to navigate from the sink (bottom right) to the saddle (center). The insets depict the corresponding vector fields in local coordinates at different iterations and also the GAD curve (orange). See also Figures~\ref{fig:isd-diffusion-maps-1}, \ref{fig:isd-diffusion-maps-5}, and~\ref{fig:isd-diffusion-maps-10}.}
  \label{fig:isd-simple-sphere-gad}
\end{figure}

The preceding algorithm applied to the vector field on the sphere $\S^2$ introduced in Example~\ref{ex:isd-simple-stereographic} produces the iterations shown in Figures~\ref{fig:isd-diffusion-maps-1}, \ref{fig:isd-diffusion-maps-5}, and~\ref{fig:isd-diffusion-maps-10}.
These figures depict the integral curves (highlighted) in the local neighborhoods obtained by sampling and integrating $\hat{X}$.
The full integral curve joining the initial point to a saddle point at the equator of the sphere in 
Figure~\ref{fig:isd-simple-sphere-gad}.

We sampled $10^3$ points per iteration and integrated the ISD vector field using an explicit Euler integrator for a total $10^3$ steps with a time-step length of $10^{-4}$.
The algorithm converges to a saddle point in ten iterations.

\subsubsection{Regular surface}

\begin{figure}[ht]
  \centering
  \includegraphics[width=\columnwidth]{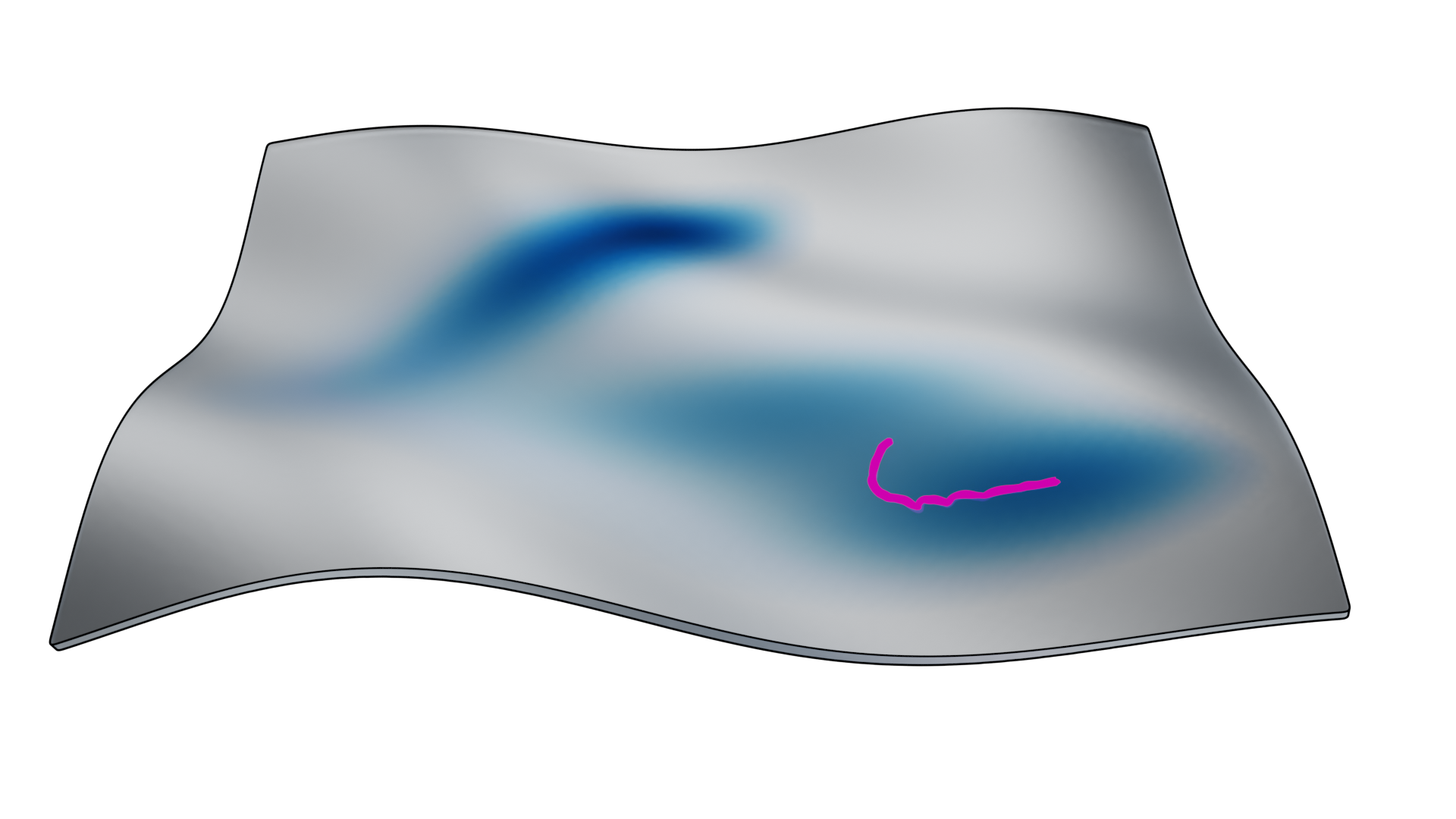}
  \caption{Path (in pink) from sink (right endpoint) to saddle point (upper endpoint) obtained by eight iterations of GAD on the M\"uller-Brown potential (blue heat map) defined on a regular surface.}
  \label{fig:muller-regular-surface}
\end{figure}

Our second example is the M\"uller-Brown (MB) potential~\cite{muller1979} mapped onto a regular surface.
Namely, the manifold $M$ is the graph of the function (see Figure~\ref{fig:muller-regular-surface}),
\begin{equation*}
  f(x^1, x^2)
  =
  \sum_{k^1 = 0}^K
  \sum_{k^2 = 0}^K
  a_{k^1, k^2}
  \cos(k^1 x^1 + k^2 x^2 + b_{k^1, k^2}),
\end{equation*}
where $K = 3$ and the coefficients are given in Table~\ref{tab:regular-surface}.
\begin{table}[ht]
  \centering
  \begin{tabular}[t]{l|l|l|l}
    $k^1$ & $k^2$ & $a_{k^1, k^2}$ & $b_{k^1, k^2}$ \\
    \hline
    \hline
    0 & 1 & 0.9490 & 0.8838 \\
    0 & 2 & 0.4575 & 0.6564 \\
    1 & 0 & 0.4152 & 0.7449 \\
    1 & 2 & 0.2911 & 0.3619 \\
    2 & 0 & 0.4121 & 0.5469 \\
    3 & 2 & 0.2817 & 0.4719 \\
  \end{tabular}
  \caption{Non-zero coefficients characterizing the regular surface.}
  \label{tab:regular-surface}
\end{table}

\begin{figure}[ht]
  \centering
  \includegraphics{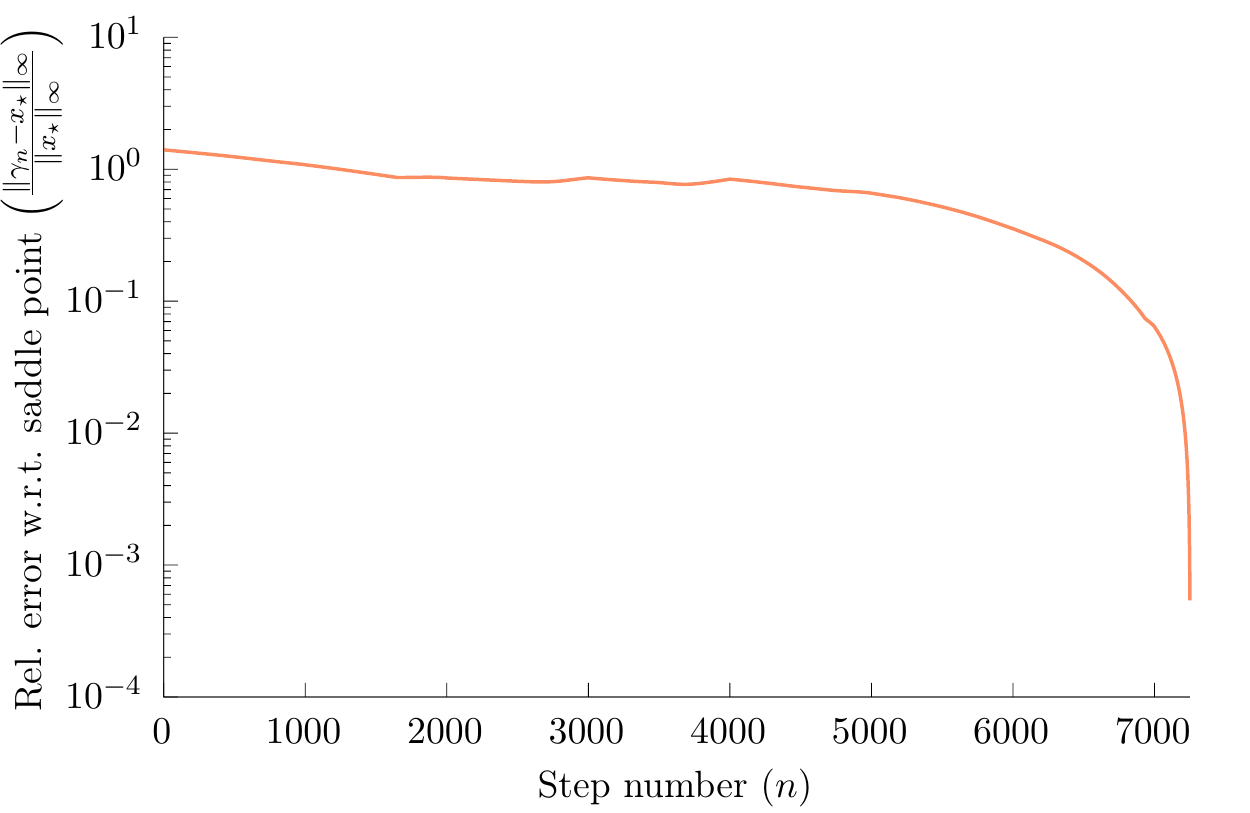}
  \caption{Relative error between each point $\gamma_n$ in the GAD path and the saddle point $x_\star$. Each iteration comprises 1000 points.}
  \label{fig:muller-path-error}
\end{figure}

After eight iterations, sampling $5 \times 10^3$ points and integrating for $10^3$ steps per iteration with a time step length of $10^{-4}$, our algorithm successfully constructs a path joining the initial point located at a minimum (sink) of the MB potential to the nearest saddle point.
The relative errors between the points in the constructed path and the known coordinates of the saddle point are shown in Figure~\ref{fig:muller-path-error}.



\section{Conclusions}
\label{sec:conclusions}

We have presented a formulation of GAD on manifolds defined by point-clouds that are meant to be sampled on-demand.
Our formulation is intrinsic and does not require the specification of the manifolds either by a given atlas or by the zeros of a smooth map.
The required charts are discovered through a data-driven, iterative process that only requires knowledge of an initial conformation of the reactant.
We illustrated our approach with two simple examples and we expect the results to transfer to the high-dimensional dynamical systems of interest in computational statistical mechanics.


\begin{acknowledgement}
  This work was supported by the US Air Force Office of Scientific Research (AFOSR) and the US Department of Energy DOE with IIT: SA22-0052-S001 and AFOSR-MURI: FA9550-21-1-0317.
\end{acknowledgement}

\appendix
\section{Appendix}

In this section we review basic notions of Riemannian geometry relevant to this paper.
Our aim is to emphasize intuition over formalism as much as possible.
We refer the reader to general treatises on the topic such as~\cite{docarmo1992,frankel2011} for a deeper presentation and we especially recommend the relevant material in~\cite{mezey1987,wales2001} for the working physical chemist.

\subsection{Tensor spaces}

Let $V$ and $W$ be two vector spaces over the real numbers and denote by $\mathscr{V}(V \times W)$ the vector space generated by all finite linear combinations of elements of the Cartesian product $V \times W$.
The \textbf{tensor product} $V \otimes W$ is a vector subspace of $\mathscr{V}(V \times W)$ with elements of the form $v \otimes w$, where $v \in V$, $w \in W$, such that:
\begin{enumerate}
\item $(v_1 + v_2) \otimes w = v_1 \otimes w + v_2 \otimes w$, where $v_1, v_2 \in V$.
\item $v \otimes (w_1 + w_2) = v \otimes w_1 + v \otimes w_2$, where $w_1, w_2 \in W$.
\item $(\alpha v) \otimes w = \alpha \, v \otimes w$, where $\alpha \in \R$.
\item $v \otimes (\alpha w) = \alpha \, v \otimes w$.
\end{enumerate}

The above construction can be recursively extended to tensor spaces of arbitrary numbers of factors (\emph{i.e.,} given vector spaces $V_1, \dotsc, V_n$, with $n \in \mathbb{N}$, we construct $V_1 \otimes V_2 \otimes \dotsb \otimes V_n$ as $V_1 \otimes (V_2 \otimes \dotsb \otimes V_n)$ and so on and so forth).

Given a vector space $V$, its \textbf{dual space}, denoted by $V^\star$, is the vector space formed by all linear maps $f \colon V \to \R$.
A tensor space of type $(p, q)$ is a tensor product
\begin{equation}
  \label{eq:tensor-space}
  T_q^p(V) = \left( \bigotimes_{i = 1}^p V \right) \otimes \left( \bigotimes_{i = 1}^q V^\star \right).
\end{equation}
Elements of $T_q^p(V)$ are called \textbf{tensors} of type $(p, q)$.

Similarly to how all vector spaces of fixed dimension over the real numbers are linearly isomorphic to each other, all tensor spaces of type $(p, q)$ over the reals are linearly isomorphic to each other and, therefore, we can assume that the factors are arranged as in~\eqref{eq:tensor-space}.

\begin{example}
  The set of linear maps from $V = \R^n$ to $W = \R^m$ is a vector space.
  If $\{ e_1, \dotsc, e_n \}$ is an orthonormal basis of $V$ and $\{ f_1, \dotsc, f_m \}$ is an orthonormal basis of $W$, then $\{ e^i \colon V \to \R \mid i = 1, \dotsc, n \}$ is the basis of the dual space $V^\star$ determined by
  \begin{equation*}
    e^i(e_j)
    =
    \begin{cases}
      1, & \text{if $i = j$,} \\
      0, & \text{otherwise}.
    \end{cases}
  \end{equation*}
  Consequently, the linear map $A \colon V \to W$ can be written as
  \begin{equation*}
    A
    =
    \sum_{i = 1}^n \sum_{j = 1}^m a_i^j \, f_j \otimes e^i,
  \end{equation*}
  where $a_i^j \in \R$.
  Additionally, if $\{ e_1, \dotsc e_n \}$ and $\{ f_1, \dotsc, f_m \}$ are, respectively, the canonical bases of $V = \R^n$ and $W = \R^m$, then $f_j \otimes e^i \in W \otimes V^\star$ coincides with the outer product $f_j e_i^\top \in \R^{m \times n}$ and $A$ coincides with the $m \times n$ real matrix
  \begin{equation*}
    A
    =
    \begin{bmatrix}
      a_1^1 & \cdots & a_1^n \\
      \vdots & \ddots & \vdots \\
      a_m^1 & \cdots & a_m^n
    \end{bmatrix}.
  \end{equation*}
\end{example}

The wedge product is defined as
\begin{equation*}
  v_1 \wedge \dotsb \wedge v_k
  =
  \frac{1}{k!}
  \sum_{\sigma \in S_k}
  \sgn(\sigma) \,
  v_{\sigma^{-1}(1)} \otimes \dotsb \otimes v_{\sigma^{-1}(k)},
\end{equation*}
where $S_k$ is the group of permutations of $\{ 1, \dotsc, k \}$ and $\sgn(\sigma)$ is equal to $-1$ if the permutation $\sigma$ is odd and $+1$ if it is even.
This is equivalent to the projection of the tensor $v_1 \otimes \dotsb \otimes v_k$ onto the subspace $\Lambda^k V$ of anti-symmetric tensors of $\bigotimes_{i = 1}^k V$.

\begin{remark}
  The signature $\sgn(\sigma)$ coincides with the Levi-Civita symbol.
  The latter should not be confused with the Levi-Civita connection, to be introduced later.
\end{remark}

\subsection{Smooth manifolds and smooth curves}

In what follows, we consider a \textbf{smooth map} to be a map that has sufficiently many derivatives and whose derivatives are continuous functions.

\begin{figure}[ht]
  \centering
  \begin{overpic}[width=0.7\columnwidth]{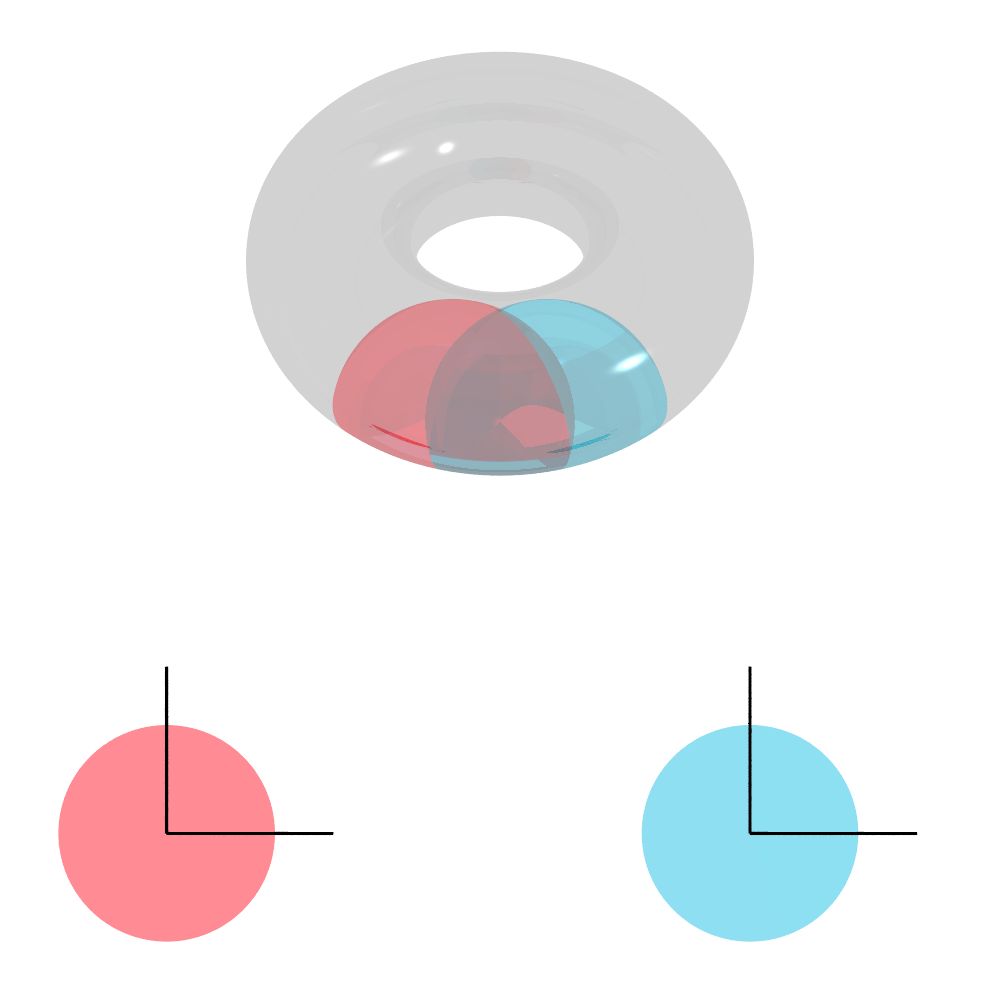}
    \thicklines
    \put(45, 45){\vector(-1, -1){15}}
    \put(55, 45){\vector(1, -1){15}}
    \put(40, 20){\vector(1, 0){20}}
    \put(32.5, 40){$\chi$}
    \put(65, 40){$\phi$}
    \put(42.5, 12.5){$\phi \circ \chi^{-1}$}
    \put(5, 30){$\R^m$}
    \put(85, 30){$\R^m$}
  \end{overpic}
  \caption{Smooth $2$-dimensional manifold (a torus) with two overlapping charts $U$ and $V$ and their corresponding systems of coordinates $\phi(U)$ and $\chi(V)$.}
  \label{fig:smooth-manifold}
\end{figure}

A \textbf{chart} is a tuple $(U, \phi)$ such that $U$ is an open set in a topological space and $\phi \colon U \to \mathbb{R}^m$ is a \textbf{diffeomorphism} (\emph{i.e.,} a smooth mapping with a smooth inverse $\psi^{-1}$).
A \textbf{smooth manifold} $M$ is a topological space together with a family of charts (\emph{i.e.,} an \textbf{atlas}),
\begin{equation*}
  \left\{ (U_\alpha, \phi_\alpha) \, | \, \phi_\alpha \colon U_\alpha \subset M \to \R^m \right\},
\end{equation*}
such that for every $p \in M$ there exists a chart $(U, \phi)$ with $p \in U $.
Moreover, for any other chart $(V, \chi)$ such that $p \in V$, the function $\phi \circ \chi^{-1}$ is smooth (see Figure~\ref{fig:smooth-manifold}).

\begin{example}
  The sets
  \begin{equation*}
    N = \left\{ (x^1, x^2, x^3) \in \R^3 \, | \, x^3 = f(x^1, x^2) = \sqrt{1 - (x^1)^2 - (x^2)^2} \right\}
  \end{equation*}
  and
  \begin{equation*}
    \mathbb{S}^2 = \left\{ (x^1, x^2, x^3) \in \R^3 \, | \, F(x^1, x^2, x^3) = 0 \right\},
  \end{equation*}
  where $F(x^1, x^2, x^3) = (x^1)^2 + (x^2)^2 + (x^3)^2 - 1$, are both smooth manifolds.
  The former is the graph of the function $f$ and the latter is implicitly defined by the zeros of the equation of a 2-sphere (\emph{i.e.,} the points in $\R^3$ at which the function $F$ vanishes).
\end{example}

Differential Geometry is often explained following the intrinsic point of view in which nothing outside of the manifold is considered.
Nevertheless, for didactic purposes it is sometimes useful to think of a manifold as being embedded in a sufficiently high-dimensional Euclidean space.
This can always be done due to the celebrated embedding theorem of Whitney~\cite{whitney1944}, which guarantees that we can embed any smooth $m$-dimensional manifold $M$ in $\R^{2 m}$.

A \textbf{smooth curve} on a manifold is a smooth mapping $\gamma \colon I \subset \R \to M$ where $I$ is an interval of $\R$.

\subsection{Tangent spaces and tangent bundles}

Again, let $M$ be a smooth $m$-dimensional manifold and let $p \in M$ be some arbitrary point on it.
A \textbf{tangent vector} $v$ at $p$ is a map that assigns a real number to each smooth, real-valued function of $M$ such that
\begin{enumerate}
\item $v(f + g) = v(f) + v(g)$,
\item $v(\lambda f) = \lambda v(f)$,
\item $v(f g) = v(f) \, g + f \, v(g)$,
\end{enumerate}
For any $f, g \colon M \to \R$ smooth and $\lambda \in \R$.

The set of tangent vectors at $p \in M$ is a vector space called the \textbf{tangent space}, denoted by $T_p M$.
Indeed, $T_p M$ is a vector space because if $v, w$ are tangent vectors at $p$ and $\lambda \in \R$, then the axioms
\begin{enumerate}
\item $(v + w)(f) = v(f) + w(f)$,
\item $(\lambda \, v)(f) = \lambda \, v(f)$,
\end{enumerate}
are satisfied.

In a local chart, the vectors
\begin{equation*}
  \frac{\partial}{\partial x^1} \bigg|_{p}, \dotsc, \frac{\partial}{\partial x^m} \bigg|_{p}
\end{equation*}
constitute a basis of $T_p M$.
\begin{example}
  Let $M = \R^m$.
  In this case, the directional derivative of $f \colon M \to \R$ at $p \in M$ in the direction $v = (v^1, \dotsc, v^m) \in \R^m$ is
  \begin{equation*}
    \nabla f(p) \cdot v = \sum_i v^i \frac{\partial f}{\partial x^i}(p).
  \end{equation*}
  In our notation, we would instead say that
  $v = \sum_i v^i \, \frac{\partial}{\partial x^i} \big|_p$
  is a tangent vector of $M$ at $p$.
\end{example}

\begin{figure}[ht]
  \centering
  \begin{overpic}[width=0.65\columnwidth,trim=1.5cm 3.5cm 1.5cm 6.5cm,clip]{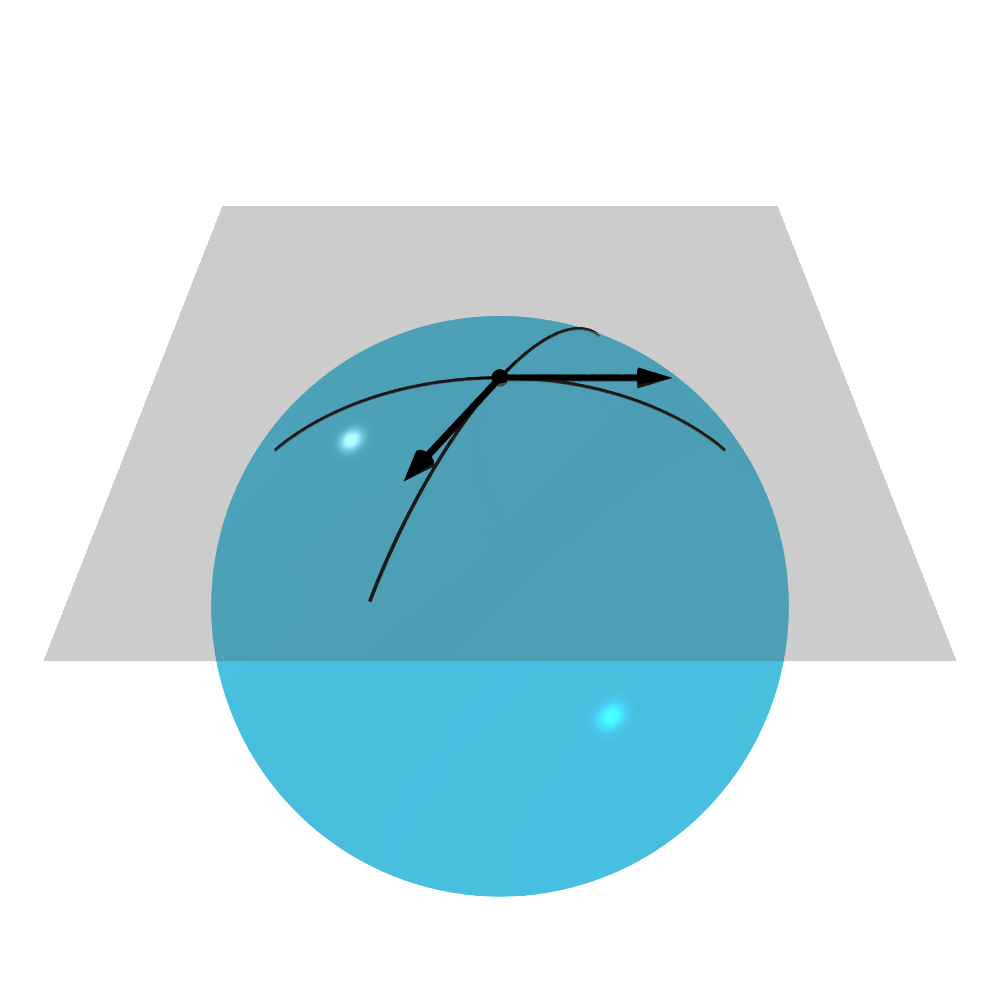}
    \put(49, 52){$p$}
    \put(83, 30){$T_pM$}
    \put(48, 10){$M$}
  \end{overpic}
  \caption{A tangent vector can be regarded intuitively as an ``arrow'' based at a point $p \in M$.
    We use tangent vectors as a tool to evaluate directional derivatives of real-valued functions defined on $M$.
    The tangent space of a manifold at $p$ can be viewed as the span of all the tangent vectors at $p$.}
  \label{fig:tangent-space}
\end{figure}

We can view each tangent vector at $p$ as the velocity vector of a smooth curve $\gamma \colon [0, 1] \subset \R \to M$ such that $\gamma(0) = p$.
In that case,
\begin{equation*}
  v(f)
  = \frac{\d}{\d t} f(\gamma(t)) = \sum_i \dot{\gamma}^i(t) \, \frac{\partial f}{\partial x^i}(\gamma(t))
  = \left( \sum_i \underbrace{\dot{\gamma}^i(t)}_{= v^i} \, \frac{\partial}{\partial x^i}\bigg|_{\gamma(t)} \right) (f).
\end{equation*}
The geometric interpretation of tangent vectors and tangent spaces is shown in Figure~\ref{fig:tangent-space}.

The collection
\begin{equation*}
  TM = \bigcup_{p \in M} T_p M
\end{equation*}
of the tangent spaces at each point of $M$ is called the \textbf{tangent bundle} and it is a smooth manifold in its own right.

\begin{example}[The tangent bundle in classical mechanics]
  If we denote the position of a mechanical system by $q \in M$, where $M$ is the configurational manifold of the system, and we consider an arbitrary trajectory $t \in [0, 1] \mapsto \gamma(t) \in M$ starting at the point $\gamma(0) = q$ with initial velocity $\dot{\gamma}(0) = v$, then we see that the tangent bundle is the collection of all the positions and velocities of the mechanical system under consideration.
  Moreover, the Lagrangian is a real-valued function of the tangent bundle,
  \begin{align*}
    L &\colon TM \to \R \\
    &(q, v) \mapsto L(q, v) = \tfrac{1}{2} v \cdot m v - U(q),
  \end{align*}
  where $m$ is the mass matrix of the system.

  The momenta $p$ are obtained from the velocities $v$ by duality (via the Legendre transform~\cite{arnold1989}).
  Because of that, the phase space of a mechanical system is a manifold that is closely related to the tangent bundle of the configurational manifold $M$.
  Indeed, the phase space is the co-tangent bundle $T M^\star$ of the configurational manifold $M$.
\end{example}

\begin{figure}[ht]
  \centering
  \begin{overpic}[width=0.5\columnwidth,trim=1.5cm 0cm 2cm 3cm,clip]{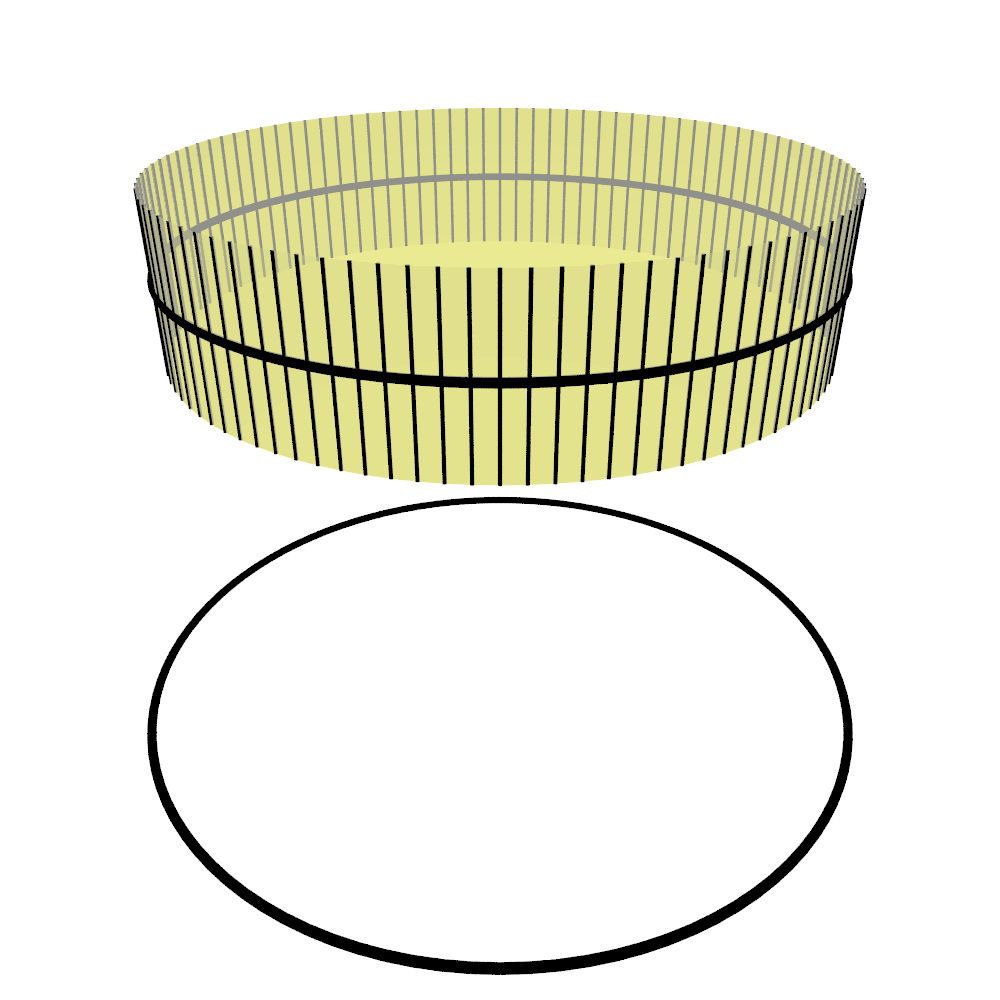}
    \thicklines
    \put(50, 60){\vector(0, -1){30}}
    \put(55, 42.5){$\pi$}
    \put(80, 5){$M$}
    \put(85, 55){$E$}
  \end{overpic}
  \caption{An example of a vector bundle.
    The base space in this case is the circle $M$ and the total space is the cylinder $E$ shown on top of $M$.
    The vertical lines on the cylinder represent the fibers (\emph{i.e.,} the one-dimensional vector spaces that are the preimages by $\pi$ of each point of $M$).}
  \label{fig:vector-bundle}
\end{figure}

The tangent bundle of a manifold is the prototype of a more general construction called a \textbf{vector bundle} that we will introduce next (albeit not fully rigorously).
Informally, a vector bundle is a collection of vector spaces that are in correspondence to each point of $M$.
Slightly more formally, we say that a \textbf{vector bundle} over $M$ is a tuple $(E, \pi, M)$ where:
\begin{enumerate}
\item The \textbf{total space} $E$ and the \textbf{base space} $M$ are smooth manifolds.
\item The \textbf{projection} $\pi \colon {E} \to {M}$ is a smooth map whose preimages (known as \textbf{fibers}) $E_p = \pi^{-1}(p)$ for $p \in M$ are vector spaces (there are additional conditions that must be satisfied but  they are beyond the scope of this appendix).
\end{enumerate}
A trivial example of a vector bundle is shown in Figure~\ref{fig:vector-bundle}.
In the case of the tangent bundle, the fiber $E_p$ is the tangent space $T_p M$ of $M$ at $p$.
More generally, we can take the fibers to be tensor spaces.

\subsection{Vector fields and sections}

\begin{figure}[ht]
  \centering
  \includegraphics[width=0.65\columnwidth,trim=1.5cm 5.5cm 1.5cm 6cm,clip]{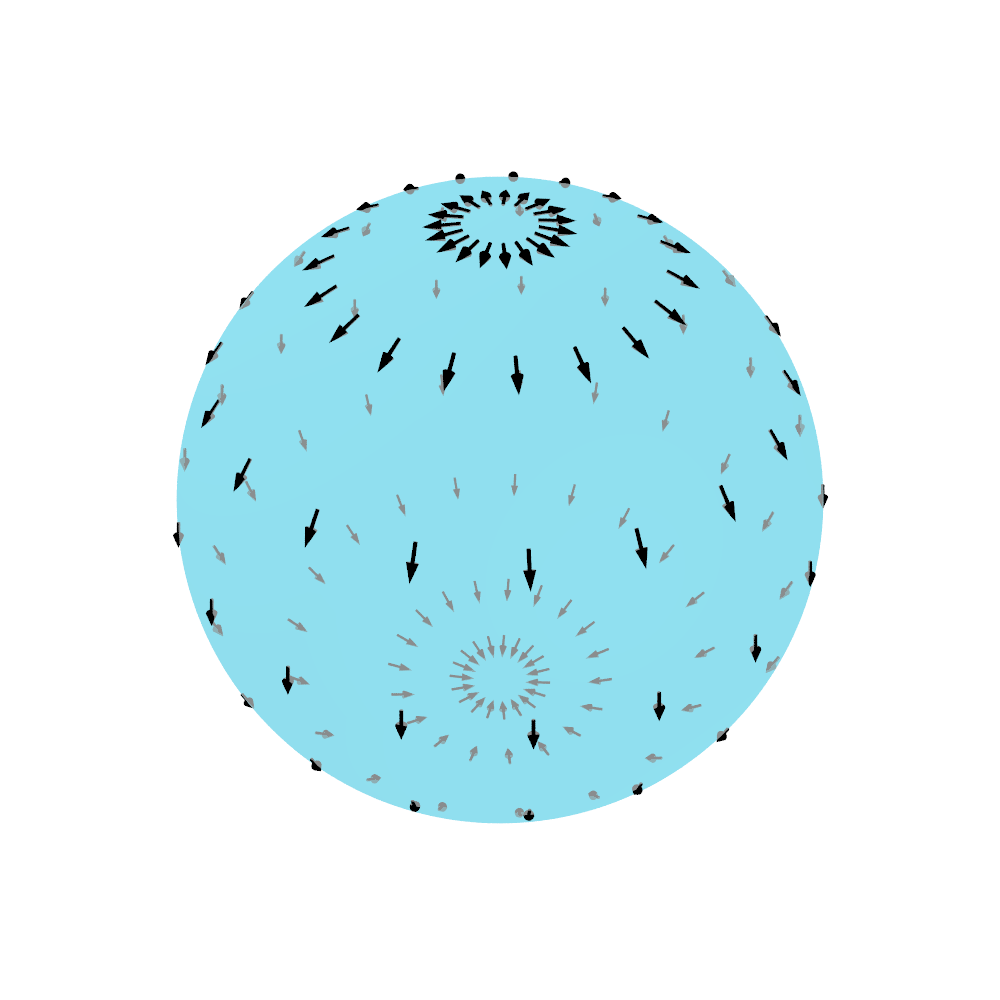}
  \caption{A vector field on the 2-sphere $\mathbb{S}^2$.
    A vector field is a (smooth) correspondence between each point $p$ of the manifold $M$ (the base of a black arrow) and a vector of $T_p M$ (the black arrow itself).}
  \label{fig:vector-field}
\end{figure}

A \textbf{vector field} is a correspondence of a tangent vector $v_p$ to each point $p \in M$.
The set of all vector fields over a manifold $M$ is denoted by $\mathfrak{X}(M)$.
An example of a vector field is shown in Figure~\ref{fig:vector-field}.

In local coordinates,
\begin{equation*}
  v_p = \sum_i v^i(p) \, \frac{\partial}{\partial x^i} \bigg|_p
\end{equation*}
From now on, we adopt the convention
\begin{equation*}
  \partial_i = \frac{\partial}{\partial x^i}.
\end{equation*}

To any smooth vector field $v \in \X(M)$ and any $t \in \R$, we can associate the \textbf{flow map}, $\Phi_t \colon M \to M$, that satisfies
\begin{equation*}
  \frac{\d}{\d t} \Phi_t(p) = v_{\Phi_t(p)}.
\end{equation*}
The smooth curve $t \mapsto \Phi_t(p)$ is also called an \textbf{integral curve} of $v$.

In local coordinates, if $v = \sum_i v^i \, \partial_i$, then $t \mapsto \Phi_t(p) = (\gamma^1(t), \dotsc, \gamma^m(t))$ solves the initial value problem
\begin{equation*}
  \left\{
    \begin{aligned}
      &\tfrac{\d}{\d t}{\gamma}^i(t) = v^i(\gamma(t)), \\
      &\text{for each $i = 1, \dotsc, m$},
    \end{aligned}
  \right.
\end{equation*}
with the initial condition $\Phi_0(p) = \gamma(0) = p$.
We often denote $\frac{\d \Phi}{\d t}$ by $\dot{\Phi}$.

Let ${\pi} \colon {E} \to {M}$ be a vector bundle.
A \textbf{section} is a smooth map
\begin{equation*}
  \begin{aligned}
    {s} \colon {M} &\to {E} \\
    {p} &\mapsto {s(p) \in E_p}
  \end{aligned}
\end{equation*}
such that the composition of $s$ followed by the projection $\pi$ is the identity in $M$ (or, equivalently, $\pi \circ s = 1_M$).
Indeed, we have the following diagram:
\begin{equation*}
  \xymatrix{
    \ar@<1ex> [d]^\pi E \\
    M \ar@<1ex> [u]^s
  }
\end{equation*}

The set of sections of a vector bundle $E$, denoted by $\Gamma(E)$, is a \emph{vector space} with the operations of addition,
\begin{equation*}
  (s + s^\prime)(p) = s(p) + s^\prime(p)
\end{equation*}
and scalar product,
\begin{equation*}
  (\lambda \, s)(p) = \lambda \, s(p)
\end{equation*}

In local coordinates $(x^1, \dotsc, x^m)$ around a point $p$ with a basis $e_1, \dotsc, e_m \in E_p$, the expression of a section $s$ can be written as
\begin{equation*}
  s = \sum_i s^i e_i.
\end{equation*}

\subsection{Pullbacks and pushforwards}

Consider a smooth map $\phi \colon M \to N$ between manifolds $M$ and $N$ with respective local coordinates $(x^1, \dotsc, x^m)$ and $(y^1, \dotsc, y^n)$.
Let $j_1, \dotsc, j_k \in \{ 1, \dotsc, n \}$ and let $f \colon N \to \R$ be a smooth function.
The \textbf{pullback} of the expression $f(y^1, \dotsc, y^n) \, \d y^{j_1} \otimes \dotsb \otimes \d y^{j_k}$ is defined by
\begin{multline*}
  \phi^\star \left(
    f(y^1, \dotsc, y^n) \, \d y^{j_1} \otimes \dotsb \otimes \d y^{j_k}
  \right) \\
  =
  f(\phi^1(x^1, \dotsc, x^m), \dotsc, \phi^n(x^1, \dotsc, x^m))
  \sum_{i_1, \dotsc, i_k = 1}^m \frac{\partial \phi^{j_1}}{\partial x^{i_1}} \dotsb \frac{\partial \phi^{j_k}}{\partial x^{i_k}} \, \d x^{i_1} \otimes \dotsb \otimes \d x^{i_k}.
\end{multline*}

Consider a vector field $X \in \X(M)$, the \textbf{pushforward} of $X$ by $\phi$ is a vector field in $\X(N)$ defined by
\begin{equation*}
  (\phi_\star X_p)(h)
  =
  X_p (h \circ \phi)
\end{equation*}
for all smooth functions $h \in C^\infty(N)$ and all points $p \in M$.

\begin{remark}
  The pullback generalizes the transformation of the integrand in the change of variables formula known from integral calculus.
  The pushforward of a vector field is an expression of the chain rule from differential calculus.
\end{remark}

\subsection{Riemannian metrics}
\label{sec:riemannian-metrics}

Let $M$ be an $m$-dimensional smooth manifold.
A \textbf{Riemannian metric} is a smooth map $g$ such that each point $p \in M$ is mapped to a symmetric and positive definite matrix $g(p) \colon T_p M \times T_p M \to \R$.
If $\partial_i, \partial_j \in T_p M$ are vectors from an orthonormal frame in the local chart $U$, the components of the Riemannian metric can be written as
\begin{equation*}
  g_{ij}(p) = g_p(\partial_i, \partial_j) = \partial_i \cdot g(p) \, \partial_j = \langle \partial_i, \partial_j \rangle,
\end{equation*}
for all $p \in U$ and $i, j = 1, \dotsc, m$.
The inner product $g_p(v, w)$ between $v = \sum_i v^i \partial_i \in T_p M$ and $w = \sum_i w^i \partial_i \in T_p M$ is defined by linearity.
The \textbf{norm} of an arbitrary vector $v = \sum_i v^i \, \partial_i \in T_p M$ is $\| v \| = \sqrt{g_p(v, v)}$.

\begin{remark}
  If we take $V = T_p M$, then its dual is $V^\star = T_p M^\star$ with orthonormal basis given by $\{ e^i \colon V \to \R \}$.
  When $V = T_p M$, as in this case, it is customary to write $\partial_i = e_i$ and $\d x^i = e^i$.
  From the previous discussion, we see that the Riemannian metric $g$ is a section in $\Gamma(T M^\star \otimes T M^\star)$ because $g(p) \colon T_p M \times T_p M \to \R$ is linear and we can write it as
  \begin{equation*}
    g(p) = \sum_{ij} g_{ij}(p) \, \d x^i_p \otimes \d x^j_p.
  \end{equation*}
\end{remark}

\begin{example}[Euclidean metric in $\R^n$]
  The Euclidean metric in $\R^n$ is
  \begin{equation}
    \label{eq:euclidean-metric}
    g = \sum_{i = 1}^n \d x^i \otimes \d x^i.
  \end{equation}
  Consider the polar coordinates in $\R^2$ given by
  \begin{equation*}
    \left\{
      \begin{aligned}
        x^1(r, \theta) &= r \cos \theta, \\
        x^2(r, \theta) &= r \sin \theta,
      \end{aligned}
    \right.
  \end{equation*}
  where $r \ge 0$ and $0 < \theta < 2 \pi$.
  We have
  \begin{equation*}
    \d x^1 = \cos \theta \, \d r - r \sin \theta \, \d \theta
    \quad \text{and} \quad
    \d x^2 = \sin \theta \, \d r + r \cos \theta \, \d \theta.
  \end{equation*}
  Therefore, the Euclidean metric~\eqref{eq:euclidean-metric} is pulled back as:
  \begin{equation*}
    g = \d x^1 \otimes \d x^1 + \d x^2 \otimes \d x^2
    = \d r \otimes \d r + r^2 \d \theta \otimes \d \theta.
  \end{equation*}
  While in Euclidean coordinates the components of the metric tensor are $g_{ij} = \delta_{ij}$, in polar coordinates they are $g_{11} = 1$, $g_{22} = r^2$, $g_{12} = g_{21} = 0$.
\end{example}

Let $I = [0, 1] \subset \R$ and let $\gamma \colon I \to M$ be a smooth curve such that $t \mapsto \gamma(t) = (\gamma^1(t), \dotsc, \gamma^m(t))$.
The \textbf{arc-length of a curve} $\gamma$ is
\begin{equation*}
  \ell(\gamma)
  =
  \int_0^1 \sqrt{\dot{\gamma}(t) \cdot g(\gamma(t)) \, \dot{\gamma}(t)} \, \d t
  =
  \int_0^1 \| \dot{\gamma}(t) \| \, \d t.
\end{equation*}

Given an arbitrary smooth manifold $M$, we can regard it as the configurational space of a (constrained) mechanical system.
Let $\gamma \colon [0, 1] \to M$ be a curve on $M$.
If we set the potential energy to be constant ($U(q) = 0$, for simplicity) and we think of the Riemannian metric $g$ as a (position-dependent) mass matrix, then it turns out that the action integral is
\begin{equation}
  \label{eq:action}
  \int_0^1 L(\gamma(t), \dot{\gamma}(t)) \, \d t,
\end{equation}
where the Lagrangian $L$ contains only a kinetic energy term,
\begin{equation*}
  L(\gamma(t), \dot{\gamma}(t)) = \tfrac{1}{2} \dot{\gamma}(t) \cdot g(\gamma(t)) \, \dot{\gamma}(t),
\end{equation*}
A curve $\gamma$ that locally minimizes the action $\int_0^1 L(\gamma, \dot{\gamma}) \d t$ is called a \textbf{geodesic curve} of $M$.
By Jensen's inequality~\cite{rudin1987}, curves that minimize~\eqref{eq:action} also minimize the arc-length.

\subsection{Connections}

Consider a vector bundle $E$ on the manifold $M$ and let $\Gamma(E)$ denote the sections of $E$.
We are interested in studying how a section $s \in \Gamma(E)$ changes along a given direction $v \in T_p M$.
In essence, we want to define what it means to take the limit
\begin{equation}
  \label{eq:alternative-covariant-derivative}
  \lim_{h \to 0} \frac{\tau_h^{-1} s(\gamma(h)) - s(\gamma(0))}{h},
\end{equation}
where $\tau_h \colon E_p \to E_{\gamma(h)}$ is a linear isomorphism between fibers and $\gamma \colon [0, 1] \to M$ is a smooth curve on $M$ such that $\gamma(0) = p$ and $\dot{\gamma}(0) = v$.
If $M = \R^n$, then we can take $\tau_h$ to be the identity map but in general we need a non-trivial $\tau_h$ to account for the fact that the fibers $E_{\gamma(h)}$ and $E_{\gamma(0)}$ are different vector spaces.
To formalize this notion, we define the concept of a connection.

Let $v, w \in \mathfrak{X}(M)$, $f \colon M \to \R$ smooth, and $\lambda \in \mathbb{R}$.
A \textbf{connection} $\nabla \colon \mathfrak{X}(M) \times \Gamma(E) \to \Gamma(E)$ acts on sections $s, t \in \Gamma(E)$ as follows:
\begin{enumerate}
\item $\nabla_v (\lambda s) = \lambda \nabla_v s$.
\item $\nabla_v (s + t) = \nabla_v s + \nabla_v t$.
\item $\nabla_v (f \, s) = v(f) \, s + f \, \nabla_v s$.
\item $\nabla_{v + w} s = \nabla_{v} s + \nabla_{w} s$.
\item $\nabla_{f v} s = f \nabla_v s$.
\end{enumerate}
In short, the connection $\nabla$ is linear on its first argument (a vector field) and acts as a derivation (\emph{i.e.,} as Leibniz's rule) on its second argument (a section).
The \textbf{covariant derivative} of the section $s$ in the direction of $v$ is defined by $\nabla_v s$.

We can characterize the covariant derivative by writing it as:
\begin{equation*}
  \nabla_v s = \bar{\nabla}_v s + A_v s,
\end{equation*}
where $\bar{\nabla}_v$ is known as the \emph{flat} connection.
If we introduce local coordinates $(x^1, \dotsc, x^m)$ around a point in the open set $U \subset M$, and consider the coordinate vector fields $\partial_k$ and a basis $e_i$ of $\Gamma(E)$, then $\bar{\nabla}_k e_i = 0$ where $\nabla_k = \nabla_{\partial_k}$.
Therefore,
\begin{equation*}
  \nabla_k e_i = \sum_j A^j_{k i} e_j,
\end{equation*}
where $A^j_{k i}$ are the components of the \textbf{vector potential}.
In full generality, we have
\begin{align*}
  \nabla_v s
  &= \nabla_{\sum_k v^k \partial_k} \left( \sum_i s^i e_i \right)
  = \sum_k v^k \nabla_{k} \left( \sum_i s^i e_i \right) \\
  &= \sum_k \sum_i v^k \left( \left( \nabla_{k} s^i \right) e_i + s^i \nabla_{k} e_i \right) \\
  &= \sum_k \sum_i v^k \left( \left( \partial_{k} s^i \right) e_i + s^i A^j_{k i} e_j \right) \\
  &= \sum_i \sum_k v^k \left( \partial_{k} s^i + s^j A^i_{k j} \right) e_i.
\end{align*}
Consequently, each component of the covariant derivative of the section $s$ along the vector $v$ is of the form
\begin{equation*}
  (\nabla_v s)^i = \sum_k v^k \left( \partial_{k} s^i + s^j A^i_{k j} \right).
\end{equation*}

\begin{remark}
  The definition of a covariant derivative $\nabla_v s$ is equivalent to~\eqref{eq:alternative-covariant-derivative} (see~\cite[Chapter 6]{spivak1999}).
  In other words, the flat connection and the vector potential fully characterize all the possible ways in which we could define a derivative as a limit of the form~\eqref{eq:alternative-covariant-derivative}.
  We will see next that there is a special type of connection that is determined by the Riemannian metric.
\end{remark}

\subsection{The Levi-Civita connection}
\label{sec:levi-civita}

A Riemannian metric induces a unique connection called the Levi-Civita connection.
The components of the vector potential of the Levi-Civita connection are called the Christoffel symbols and we shall derive them next using classical mechanics.

Consider a free particle with unit mass moving on the manifold $M$ with position $q$ and velocity $\dot{q}$.
Its Lagrangian is
\begin{equation*}
  L(q, \dot{q})
  =
  \frac{1}{2} \sum_{i = 1}^m \sum_{j = 1}^m g_{ij}(q) \, \dot{q}^i \dot{q}^j.
\end{equation*}
The Euler-Lagrange equations for this mechanical system are
\begin{equation*}
  \frac{\d}{\d t}\frac{\partial L}{\partial \dot{q}^i} - \frac{\partial L}{\partial q^i} = 0,
  \quad \text{for} \quad i = 1, \dotsc, m.
\end{equation*}
This implies that
\begin{align*}
  0 &=
      \sum_{j = 1}^m \frac{\d}{\d t} \left( g_{ij}(q) \, \dot{q}^j \right)
      -
      \frac{1}{2} \sum_{j = 1}^m \sum_{k = 1}^m \frac{\partial  g_{jk}}{\partial q^i}(q) \, \dot{q}^j \dot{q}^k \\
    &=
      \sum_{j = 1}^m g_{ij}(q) \, \ddot{q}^j
      +
      \sum_{j = 1}^m \sum_{k = 1}^m \left( \frac{\partial  g_{ij}}{\partial q^k}(q) - \frac{1}{2} \frac{\partial  g_{jk}}{\partial q^i}(q) \right) \dot{q}^j \dot{q}^k,
\end{align*}
for $i = 1, \dotsc, m$.
Using the identity
$\sum_{j = 1}^m \sum_{k = 1}^m \frac{\partial g_{ij}}{\partial q^k}
  =
  \sum_{j = 1}^m \sum_{k = 1}^m \frac{\partial g_{ik}}{\partial q^j},$
we find
\begin{equation*}
  0
  =
  \sum_{j = 1}^m g_{ij}(q) \, \ddot{q}^j
  +
  \frac{1}{2} \sum_{j = 1}^m \sum_{k = 1}^m \left( \frac{\partial  g_{ij}}{\partial q^k}(q) + \frac{\partial  g_{ik}}{\partial q^j}(q) - \frac{\partial g_{jk}}{\partial q^i}(q) \right) \dot{q}^j \dot{q}^k,
\end{equation*}
for $i = 1, \dotsc, m$.
Finally, multiplying both sides above by $\sum_{i = 1} g^{\ell i}$, where $g^{\ell i}$ are the components of the inverse metric tensor (\emph{i.e.,} $\sum_{i = 1} g^{\ell i} g_{i j} = \delta_{\ell, j}$ where $\delta_{\ell j}$ is the Kronecker delta), we arrive at the geodesic equation
\begin{equation}
  \label{eq:geodesic-equation}
  0
  =
  \ddot{q}^\ell
  +
  \sum_{j = 1}^m \sum_{k = 1}^m
  \Gamma^\ell_{jk}(q) \, \dot{q}^j \dot{q}^k,
\end{equation}
where
\begin{equation*}
  \Gamma^\ell_{jk}(q)
  =
  \frac{1}{2} \sum_{i = 1}^m g^{\ell i}(q) \left( \frac{\partial  g_{ij}}{\partial q^k}(q) + \frac{\partial  g_{ik}}{\partial q^j}(q) - \frac{\partial g_{jk}}{\partial q^i}(q) \right)
\end{equation*}
is the expression of the Christoffel symbol for $i, j, \ell \in \{ 1, \dotsc, m \}$.

\begin{remark}
  Equation~\eqref{eq:geodesic-equation} is known as the \textbf{geodesic equation} and can be written more succinctly as $\nabla_{\dot{q}} \dot{q} = 0$.
  The phase flow of this equation when the initial velocity has unit length is known as the \textbf{exponential map} $\exp_t$ for $t \ge 0$.
  Therefore, $\exp_0 p = p$, $\frac{\d}{\d t} \exp_t p = \dot{q}(t)$ such that $\| \dot{q}(0) \| = 1$, and $q(t) = \exp_t p$ solves the second order ODE~\eqref{eq:geodesic-equation}.
\end{remark}

\begin{remark}
  In the Euclidean case, $M = \R^m$ and $g_{ij} = \delta_{ij}$, leading to the usual expression of the kinetic energy and to the free particle moving in a straight line.
  When $M$ is an arbitrary Riemannian manifold and the free particle moves along a geodesic, we can interpret the Christoffel symbols as the terms giving rise to the corrections to an otherwise straight trajectory so that it remains on the manifold $M$.
\end{remark}

\subsection{Raising and lowering indices}

The inner product $g \colon T M \times T M \to \mathbb{R}$ between two vectors $X, Y \in T_p M$ determined by the Riemannian metric gives rise to a linear isomorphism (called \emph{flat}) between $T_p M$ and $T_p^\star M$ by mapping
\begin{equation*}
  X \in T_p M \mapsto X^\flat = g(X, \cdot) \in T_p^\star M.
\end{equation*}
The inverse linear map (called \emph{sharp}) is
\begin{equation*}
  \omega \in T_p^\star M \mapsto \omega^\sharp \in T_p M
\end{equation*}
such that $\omega = g(\omega^\sharp, \cdot)$.
The isomorphisms introduced above are called the \emph{musical isomorphisms}.

Denoting by $(g^{ij})$ the inverse of the metric tensor $g$, we write the components of $\omega^\sharp$ and $X^\flat$ as
\begin{equation*}
  {(\omega^\sharp)^i = \sum_j g^{ij} \omega_j}
  \quad \text{and} \quad
  {(X^\flat)_i = \sum_j g_{ij} X^i}.
\end{equation*}

\begin{example}
  The gradient of a function $f \in C^\infty(M)$ is defined as $\grad f = \d f^\sharp$, where $\d f = \sum_{i = 1}^m \frac{\partial f}{\partial x^i} \d x^i$.
  We have seen that in the case of the plane $\R^2$ with polar coordinates, $g = \d r \otimes \d r + r^2 \d \theta \otimes \d \theta$.
  Consequently, $\grad f = \frac{\partial f}{\partial r} \frac{\partial}{\partial r} + \frac{1}{r^2} \frac{\partial f}{\partial \theta} \frac{\partial}{\partial \theta}$.
\end{example}

The musical isomorphisms can be used with individual factors of a tensor product to map between the primal and dual vector spaces.


\bibliography{bibliography}

\end{document}